\documentstyle[epsf]{article}

\newtheorem{example}{Example}[section]
\newtheorem{note}[example]{Note}
\newtheorem{theorem}[example]{Theorem}
\newtheorem{corollary}[example]{Corollary}
\newtheorem{definition}[example]{Definition}

\newtheorem{lemma}[example]{Lemma}

\newcommand{\bi}{\begin{itemize}}
\newcommand{\ei}{  \end{itemize}}
\newcommand{\be}{\begin{equation}}
\newcommand{\ee}{  \end{equation}}
\newcommand{\bea}{\begin{eqnarray}}
\newcommand{\eea}{  \end{eqnarray}}

\def\wt{{\rm wt\,}}
\def\Proof{\medskip\noindent {\it Proof: }}
\def\cf#1#2{(#1\!:\!#2)}
\def\cqfd{\hfill $\Box$ \medskip}
\def\endex{\hfill $\diamond$ \medskip}

\def\mathalign#1{\hbox to 0pt{\hss$\vcenter{\openup\jot
\tabskip=0pt plus1fil \halign to \hsize
{\hfil$\displaystyle ##$\tabskip=0pt&$\displaystyle ##$\hfil
&&$\displaystyle ##$\hfil
\tabskip=0pt plus1fil \cr #1}}$\hss}}

\def\makestrut#1#2{{\dimen12=#2
\divide\dimen12 by 4\dimen11=\dimen12\multiply\dimen11 by 3
\global\setbox#1=\hbox{\vrule height\dimen11 depth\dimen12 width0pt}}}

\newdimen\tadhdimen \newdimen\tabhdimen \newdimen\vdimen
\newdimen\smtadhdimen \newdimen\smtabhdimen
\newcount\splurt
\newbox\tadstrut \newbox\tabstrut
\newbox\smtadstrut \newbox\smtabstrut

\def\setyoungsize#1#2{          
  \tadhdimen=#1\tabhdimen=#1\advance\tabhdimen by -0.4truept%
  \vdimen=#2%
  \makestrut\tadstrut\vdimen
  \advance\vdimen by -0.4pt%
  \makestrut\tabstrut\vdimen}

\def\setsmyoungsize#1#2{        
  \smtadhdimen=#1\smtabhdimen=#1\advance\smtabhdimen by -0.4truept%
  \vdimen=#2%
  \makestrut\smtadstrut\vdimen
  \advance\vdimen by -0.4pt%
  \makestrut\smtabstrut\vdimen}

\setyoungsize{15pt}{15pt}
\setsmyoungsize{9pt}{9pt}

\def\youngt#1{%
  \vcenter{\offinterlineskip
  \halign{&\copy\tadstrut\hbox to \tadhdimen{\hss$##$\hss}\cr #1}}}
\def\youngd#1{%
  \vcenter{\offinterlineskip
  \halign{&\vrule##&\copy\tabstrut\hbox to \tabhdimen{\hss$##$\hss}\cr #1}}}

\def\smyoungt#1{{\vcenter{\offinterlineskip
  \halign{&\copy\smtadstrut
  \hbox to \smtadhdimen{\hss$\scriptstyle ##$\hss}\cr #1}}}}
\def\smyoungd#1{{\vcenter{\offinterlineskip
  \halign{&\vrule##&\copy\smtabstrut
  \hbox to \smtabhdimen{\hss$\scriptstyle ##$\hss}\cr #1}}}}

\def\hdashfill{\leaders\hbox to 3pt{%
\hfil\vrule width1.5pt height0.4pt depth0pt}\hfill}

\def\gauss#1#2{{\textstyle\left[{#1\atop #2}\right]}}
\def\qeta#1{{\textstyle{1\over(q)_{#1}}}}

\begin{document}

\title{A Burge tree of Virasoro-type polynomial identities}

\author{
Omar~{\sc Foda}\thanks{foda@maths.mu.oz.au},\  
Keith~S.~M.~{\sc Lee}\thanks{ksml@maths.mu.oz.au},\ 
and\ 
Trevor~A.~{\sc Welsh}\thanks{trevor@maths.mu.oz.au}\\
}
\date{
{\it Department of Mathematics and Statistics, \\
     The University of Melbourne, \\
     Parkville, Victoria 3052, \\
     Australia. \\}
}

\maketitle

\begin{abstract}
Using a summation formula due to Burge, and a combinatorial 
identity between partition pairs, we obtain an infinite tree 
of $q$-polynomial identities for the Virasoro characters
$\chi^{p, p'}_{r, s}$, dependent on two finite size
parameters $M$ and $N$, in the cases where:

\begin{enumerate}
\item $p$ and $p'$ are coprime integers that satisfy
      $0 < p < p'$. 
\item If the pair $\cf {p'}p$ has a continued fraction 
      $(c_1, c_2, \ldots, c_{t-1}, c_t+2)$, where $t \ge 1$,
      then the pair $\cf sr$ has a continued fraction 
      $(c_1, c_2, \ldots, c_{u-1}, d)$, 
      where $1 \le u \le t$, and $1 \le d \le c_{u}$.
\end{enumerate}

The limit $M \to \infty$, for fixed $N$, and the limit 
$N \to \infty$, for fixed $M$, lead to two independent 
boson-fermion-type $q$-polynomial identities: in one case, 
the bosonic side has a conventional dependence on 
the parameters that characterise the corresponding character. 
In the other, that dependence is not conventional. In each case, 
the fermionic side can also be cast in either of two different 
forms. 

Taking the remaining finite size parameter to infinity in either 
of the above identities, so that $M \to \infty$ and $N \to \infty$, 
leads to the same $q$-series identity for the corresponding 
character.

\end{abstract}

\newpage

\setcounter{secnumdepth}{10}
\setcounter{section}{-1}

\section{Introduction}

\subsection{A \lq lattice\rq\ of Virasoro characters}

Let us consider the set of minimal conformal field theories $M(p, p')$ 
of Belavin, Polyakov and Zamolodchikov \cite{bpz}\footnote{For an 
introduction to conformal field theory, see \cite{dms-book}.}. They 
are labelled by two coprime integers\footnote{Two integers that have 
no common divisor.} $\{p, p'\}$, where $0 < p < p'$. If we think of 
$\{p, p'\}$ as the coordinates of points on a two-dimensional 
square lattice, then crudely speaking, one can say that there is 
a two-fold infinity of these theories\footnote{Of course, this is an 
inaccurate description as one can order all sets with a finite number 
of elements that take discrete values on the line. However, this 
description will be useful.}. 
For $p' = p+1$, they 
correspond to the critical limit of the lattice models of Andrews, 
Baxter and Forrester \cite{abf}. For $p' > p+1$, they correspond to 
the critical limit of the models of Forrester and Baxter 
\cite{fb}\footnote{For an introduction to lattice models, see 
\cite{baxter-book}.}. 

The spectrum of a minimal theory $M(p, p')$ on the cylinder, or 
equivalently, the set of one-point functions of the corresponding 
lattice model on the plane, can be written in terms of Virasoro 
characters $\chi^{p,p'}_{r,s}$, where 
$1 \le r < p$, $1 \le s < p'$, and 
$\chi^{p,p'}_{r,s} = 
 \chi^{p,p'}_{p-r,p'-s}$\cite{djkmo}\footnote{To be precise, 
the one-point functions are proportional to Virasoro characters 
\cite{djkmo}.}. In the sense used above, one can say that there is 
a four-fold infinity of these characters.

\subsection{Q-series identities}

The characters $\chi^{p,p'}_{r,s}$ are $q$-series. The Stony Brook group 
were the first to realise that the form of these $q$-series is not unique
\cite{stony-brook-review}\footnote{In the context of affine algebras, a 
similar observation was made and used in \cite{lepowsky}.}.
The different forms arise naturally from
the different approaches to computing them.
Each approach arises from a specific physical
interpretation\footnote{For an introduction and a review
of the physics behind the different forms, we refer to 
\cite{stony-brook-review}.}. All we need to know here, is that depending 
on how they are computed, the characters can be
\lq alternating-sign\rq\ series, or \lq constant-sign\rq\ series.

For physical reasons \cite{stony-brook-review},
the alternating-sign expressions are called \lq bosonic expressions\rq.
Those of constant-sign are called \lq fermionic expressions\rq.
For all $p$, $p'$, $r$ and $s$, bosonic expressions for
$\chi^{p,p'}_{r,s}$ are given by Rocha-Caridi \cite{rocha}:
\begin{equation}\label{RochaEq}
\chi^{p, p'}_{r, s}=
{1\over (q)_\infty}\sum_{k=-\infty}^\infty
(q^{k^2pp'+k(p'r-ps)}-q^{(kp+r)(kp'+s)}).
\end{equation}
However, fermionic forms for $\chi^{p,p'}_{r,s}$ are known explicitly
only for certain $p$, $p'$, $r$ and $s$.
In those cases for which they are known, equating the two forms
results in $q$-series identities.
For example, a fermionic expression for $\chi^{2,5}_{1,2}$
is provided by:
\begin{equation}\label{RReq}
\chi^{2,5}_{1,2}=
\sum_{n=0}^\infty {q^{n^2}\over (q)_n}\;.
\end{equation}
Equating this expression with the relevant instance of
(\ref{RochaEq}) results in one of the celebrated
Rogers-Ramanujan (RR) identities.
By the same means, if we knew fermionic expressions for
all $\chi^{p,p'}_{r,s}$, then
we would end up with a four-fold infinity of generalised RR 
identities\footnote{The original Rogers-Ramanujan identities
involve two equalities: an equality between a constant-sign series and
an alternating-sign series, and an equality between the latter and 
a product form. Strictly speaking, it is the latter equality that
is challenging from a combinatorial point of view. We do not
consider the product form in this work.}.

\subsection{$Q$-polynomial identities}

Rather than work in terms of $q$-series, it is possible to work in 
terms of $q$-polynomials. The latter arise by taking the
size of underlying physical models to be finite and to depend
on various {\em finite size} parameters.
In the limit for which the finite size parameters tend to infinity,
we recover the original $q$-series.  

Working in terms of $q$-polynomials is particularly suited to the
combinatorial approach that we follow in this work: one considers 
a class of combinatorial objects, $q$-counts them in two different 
ways\footnote{By $q$-counting, we mean counting objects in such a 
way that one keeps track of a certain statistic that we call the 
\lq weight\rq. In the case of partitions, the weight is simply the 
integer that is partitioned.} and then equates the results. In this 
work, the combinatorial objects that we consider are pairs of partitions 
that satisfy certain conditions\cite{andrews-red-book,abbbfv,burge}. 

As mentioned above, bosonic expressions are known for all 
$\chi^{p,p'}_{r,s}$. It is also straightforward to obtain
finite analogues of $\chi^{p,p'}_{r,s}$ that depend on one finite
size parameter \cite{abbbfv}. Furthermore, finite analogues
that depend on more than one finite size parameter are known
\cite{burge}, as we will see below.  

One approach to obtaining new polynomial RR identities is to explicitly 
evaluate more generating functions in fermionic form. Another approach 
is to use summation formulae and $q$-series transforms to generate new 
expressions and identities from known ones. 

\subsection{The Bailey transform}

The Bailey transform \cite{bailey}, with extensions by Andrews 
\cite{andrews-on-bailey,andrews-white-book}, has been used to generate 
infinite sequences of new identities from known ones \cite{fq,bms}. 
An initial known identity acts as a starting point, or as a \lq seed\rq\,
for an infinite sequence of identities. The sequence obtained is 
one-fold infinite and covers a measure zero subset
of the full set of possible identities. 

That the Bailey transform allows one to obtain more \lq complicated\rq\ 
identities (in the sense of larger $p$ and $p'$) from \lq simpler\rq\ ones
(in the sense of smaller $p$ and $p'$), raises the possibility 
that one can actually obtain {\it all} identities, for all 
Virasoro characters from a single {\it combinatorially trivial} 
identity. However, it is unclear how to this may be achieved in
the context of the Bailey transform, as will be discussed below. 

\subsection{The Burge transform}

In this work, we take a step towards generating the entire set of 
RSOS character identities from a single simple one. We make use of 
a transform due to Burge which, in a sense, generalises a restricted 
version of the Bailey transform: the Bailey transform involves 
two continuous parameters and the Burge transform generalises the special 
case of the Bailey transform where both parameters tend to infinity. 

\subsection{A comparison of two transforms}

The Burge transform is \lq stronger\rq\ than the Bailey transform in the 
following sense: in both cases, one starts from $q$-polynomial identities.  
In the case of the Bailey transform, the polynomials that the transform 
acts on to generate new identities, depend on a single {\it finite 
size} parameter, say $L$\footnote{Roughly speaking, this is related 
to the fact that the objects that are counted are single partitions.}. 
The action of the Bailey transform is such that the final result is 
not a polynomial identity but a $q$-series identity.

In the case of the Burge transform, the transform acts on $q$-polynomial 
identities that depend on {\it two} finite size parameters, say $M$, 
and $N$\footnote{Roughly speaking, this is related to the fact that the 
objects that are counted are are pairs of partitions.}${}^{,}\!$
\footnote{To be more 
precise, the Burge transform involves {\it four} finite size parameters, 
say $N$, $M$, $N'$, and $M'$. However, the identities that we obtain can 
all be derived in terms of two parameters only: $M=M'$, and $N=N'$. The 
general case of four parameters is relevant to identities that correspond 
to the most general RSOS characters. We do not deal with the most general 
case in this work.} 

As noticed by Burge \cite{burge}, working in terms of partition pairs 
allows more \lq games\rq\ to be played:
it allows us to find more transformations 
under which the generating function of the objects that are counted remains 
invariant. This extra freedom is what allows us to use the Burge transform
to obtain a two-fold infinity of identities rather than a one-fold infinity 
as in the case of the Bailey transform. It is also what allows us to obtain
two independent identities for each character. 

Furthermore, the result of applying the Burge transform to a $q$-polynomial 
with two (four) finite size parameters is once again a $q$-polynomial with 
two (four) finite size parameters.

\subsection{Outline of results}

The multi-parameter polynomials involved in the Burge transform are 
not precisely finite analogues of Virasoro characters. To obtain
the latter, one has to take one of the finite parameters, say
$M$ to infinity, and normalise the resulting expressions
appropriately. Only then does one obtain a finite analogue of 
a Virasoro character. Because $p$ can be any positive integer, and 
$p'$ can be any positive integer that is larger than $p$ and coprime
to it, we obtain a two-fold infinite tree of polynomial
identities. If one takes the other finite parameter $N$ to
infinity, instead of $M$, one obtains another two-fold infinite 
tree. Taking the finite parameter that is left to infinity, one
obtains $q$-series identities for the characters.

We proceed in two steps. Firstly, we use the Burge transform, without
any further additions, to obtain an identity for {\it one} character
of each of the $M(p, p')$ models: $\chi^{p,p'}_{r_0,s_0}$, where
$| p s_0 - p' r_0 |=1$\footnote{This is the character that has the
smallest conformal dimension in the model\cite{bm}.}.

Secondly, we make use of a simple combinatorial identity, to extend
our results to further characters $\chi^{p, p'}_{r_i,s_i}$ as will
be described below. What we can say at this point is the following:
if $p'/p$ has a continued fraction
expansion of the form $(c_1, c_2, \ldots, c_{t-1}, c_t+2)$ with $t \ge 1$,
then $s_i/r_i$ has a continued fraction expansion of the form
$(c_1, c_2, \ldots, c_{u-1}, d)$, where $1 \le u \le t$ and
$1 \le d \le c_{u}$
\footnote{
We can say that the latter correspond to all
possible \lq tapered truncations\rq\ of the former.}.
(As we will see, the case $\{r_0,s_0\}$ is included in this set.)

\subsection{Remarks on content}

At this point, we wish to emphasise that, mathematically speaking, 
this work is based entirely on the ideas of Burge, as expressed in 
\cite{burge}. In particular, the idea of generating a tree of polynomial
identities that depend on two finite size parameters, and that reduce 
to two independent polynomial identities by taking one of the two 
parameters to infinity is contained explicitly in \cite{burge}. We 
believe that even the combinatorial identity that we use to obtain 
an extended Burge tree, is also known to Burge, although it was not 
explicitly stated in \cite{burge}\footnote{There are two results in 
\cite{burge} that, to the best of our understanding, could only be 
derived using such an identity.}. 

The first purpose of this work is to introduce and clarify Burge's 
ideas, which are often cited, but seldom read\footnote{We find that 
though \cite{burge} is ingenious, it is also very succinctly written, 
and therefore not easy to read. In particular, certain consistency 
conditions that must be imposed on the partitions pairs are not 
explicitly stated. This difficulty is further compounded by the fact 
that it contains, unfortunately, a large number of misprints.}.
For this reason we have included derivations of all of the relevant 
results from \cite{burge}.

Furthermore, we wish to make an explicit connection with the
physically-motivated works on boson-fermion Virasoro character
identities. In particular, we wish to use Burge's results to
obtain algorithms to generate the largest possible set of
character identities. The emphasis in this work is on 
systematic and algorithmic methods to generate characters.  

\subsection{Outline of paper}

In \S 1, we introduce the combinatorial objects that we are
interested in: partition pairs that obey specific restrictions. 
In \S 2, we derive a $q$-polynomial summation formula due to 
Burge. We refer to this formula as the \lq Burge transform\rq.
In \S 3, we explain Burge's algorithm for generating a tree of 
$q$-polynomial identities using the Burge transform. The main 
result here, Theorem \ref{DetOneThrm}, gives finitised character 
identities for the $\chi^{p,p'}_{r_0,s_0}$ mentioned above.

In \S 4, we introduce a combinatorial identity that enables the 
algorithm of the previous section to further extend the tree of 
polynomial identities. The remaining sections are devoted to our 
results. In \S 5, we obtain a tree of polynomial identities that 
depend on one finite size parameter and that have a bosonic side 
with a conventional dependence on the parameters that characterise 
the corresponding character. The resulting one-parameter finitised
$q$-polynomial Virasoro characters are unified and presented in 
Corollary \ref{FirstCor}.

In \S 6, we obtain a tree of polynomial identities that depend on 
one finite size parameter and that have a bosonic side with a 
non-conventional dependence on the parameters that characterise
the corresponding character. Here, the resulting one-parameter 
finitised $q$-polynomial Virasoro characters are presented in 
Corollary \ref{SecondCor}. In \S 7, we remove both finite size 
parameters, and obtain a tree of character identities. The 
resulting Virasoro character ($q$-series) identities are unified 
and presented in Corollary \ref{ThirdCor}. In \S 8, we include 
a number of remarks.

Appendix A details Burge's proof of his alternating sign generating 
function for the restricted partition pairs. Appendix B provides 
some details for examples given in \S 6 and \S 7.

\section{Restricted partition pairs}

\subsection{Combinatorial objects related to Virasoro characters}

Although we are ultimately interested in generating identities that 
express Virasoro characters in two ways with different
physical interpretations, the approach followed in this work is 
combinatorial in the following sense: following Burge \cite{burge}, 
we obtain a tree of $q$-polynomial identities by enumerating
certain combinatorial objects in two different ways,
and equating the results\footnote{To be more precise, 
we will enumerate these objects only in a special case where the 
conditions they satisfy are so strict that the set of enumerated 
objects can only be empty. Once we obtain an identity that expresses 
the above fact, we use it as a root of a tree of identities, and 
generate the rest of the tree using the Burge transform, and a simple 
combinatorial identity between partitions.}. 

In previous work, various combinatorial objects have been used in order 
to obtain finite versions of Virasoro characters:

\begin{enumerate}
\item{Paths} \cite{fw,w1,w2}
\item{Partitions with prescribed hook differences} \cite{fq,fw}
\item{Coloured Young diagrams} \cite{fow} 
\item{Tableaux} \cite{kkr,df} 
\end{enumerate}

In this work, we use yet another type of combinatorial object,
introduced by Burge \cite{burge}, that are most closely
related to the partitions with prescribed hook differences that
were introduced in \cite{abbbfv}.

\subsection{Definitions}

A partition $p=(p_N,p_{N-1},p_{N-2},\ldots,p_1)$ in $(N,M)$ is a
sequence of $N$ integer parts $\{p_N,p_{N-1},p_{N-2},\ldots,p_1\}$
such that $M\ge p_N\ge p_{N-1}\ge\cdots\ge p_1\ge0$\footnote{Note 
that, following Burge, the parts here are labelled unconventionally,
with decreasing labels.}.
The weight $\wt(p)$ of $p$ is the sum of its parts:
$\wt(p)=p_N+p_{N-1}+p_{N-2}+\cdots+p_1$.
Equivalently, $\wt(p)$ is the number of nodes in the corresponding
Young diagram of $p$.

A partition pair $(q,p)$ in $(N_1,M_1)\times(N_2,M_2)$ is an ordered 
pair such that $q$ is a partition that fits in a box of dimensions 
$(N_1,M_1)$ and $p$ is a partition that fits in a box of dimensions 
$(N_2,M_2)$.

A partition pair can be depicted as follows: Firstly, we draw an 
$N_2\times N_1$ 'Durfee rectangle'. Next, we attach the $q$ partition 
to the bottom edge, and attach the $p$ partition to the right edge, 
as indicated in the figure below.

\medskip
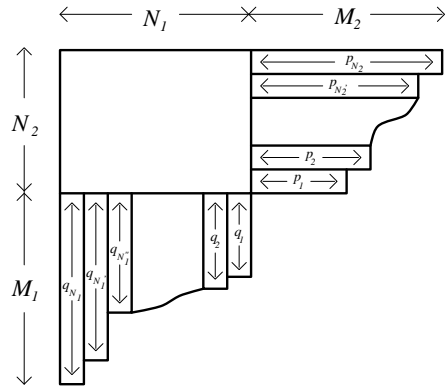
\centerline{\epsfbox{pair.eps}}
\medskip

\noindent
(Here $N^\prime=N-1$ and $N^{\prime\prime}=N-2$.)
The weight $\wt((q,p))$ of the partition pair $(q,p)$ is given by 
$\wt((q,p))=\wt(q)+\wt(p)$. In other words, the Durfee rectangle is 
empty, and does not contribute to the weight of the partition pair. 

Given a set of non-negative integer parameters 
$\{N_1,M_1,N_2,M_2,a,b,\alpha,\beta\}$ that satisfy the
conditions\footnote{The last two of the following four
conditions, and the special cases that follow, are not 
explicitly stated in \cite{burge}. However, they are required 
for consistency\cite{kg}.} 

\begin{enumerate}
\item $a+b>0$, 
\item $\alpha+\beta>0$, 
\item $-a\le N_1-N_2\le b$,
\item $-\beta\le M_1-M_2\le\alpha$,
\end{enumerate}
we say that a partition pair $(q,p)$ in $(N_1,M_1)\times(N_2,M_2)$ 
is restricted if
\begin{equation}\label{restrictions}
\mathalign{
p_i-q_{i+1-a}\: &\ge 1-\alpha,&\cr
q_i-p_{i+1-b}\: &\ge 1-\beta,&\cr}
\end{equation}
for all meaningful values of the subscripts --- that is, for
$\max\{1,a\}\le i\le\min\{N_2,N_1+a-1\}$ in the first inequality and for
$\max\{1,b\}\le i\le\min\{N_1,N_2+b-1\}$ in the second --- and the following
special cases also hold:
\begin{equation}\label{exceptions}
\mathalign{
&\hbox{if $a=0$ then }\quad &q_1\le\alpha-1;\cr
&\hbox{if $b=0$ then }\quad &p_1\le\beta-1;\cr
&\hbox{if $a=N_2-N_1$ then }\quad &p_{N_2}\ge M_1-\alpha+1;\cr
&\hbox{if $b=N_1-N_2$ then }\quad &q_{N_1}\ge M_2-\beta+1.\cr}
\end{equation}

We define
${\cal R}(N_1,M_1,N_2,M_2,a,b,\alpha,\beta)$
to be the set of all such restricted partition pairs.
The generating function for these pairs is then defined by:
$$
R(N_1,M_1,N_2,M_2,a,b,\alpha,\beta;q)=
\sum_{{\cal R}(N_1,M_1,N_2,M_2,a,b,\alpha,\beta)}
q^{\wt((q,p))}.
$$
For convenience, we can drop the last parameter in $R$, and
define 
$$
R(N_1,M_1,N_2,M_2,a,b,\alpha,\beta)=
R(N_1,M_1,N_2,M_2,a,b,\alpha,\beta;q).
$$

For the purposes of the result given below, we define the
$q$-shifted factorial $(q)_n$, and the Gaussian polynomial 
$\left[{P\atop N}\right]$ as follows
\begin{eqnarray}
(q)_n&=&(1-q)(1-q^2)(1-q^3)\cdots(1-q^n)\qquad(n\ge0),\\[0.5mm]
\left[{P\atop N}\right]&=&
\cases{
\displaystyle{(q)_P\over (q)_N (q)_{P-N}}&if $P\ge N\ge0$;\cr
\noalign{\smallskip}
0&otherwise.\cr}\label{GaussDef}
\end{eqnarray}
Further, we define
$$
g(N_1,M_1,N_2,M_2,x,y)=\left[{N_1+M_1+x-y\atop N_1+x}\right]
\left[{N_2+M_2-x+y\atop N_2-x}\right].
$$

\subsection{Bosonic generating function}

Now that we have defined the combinatorial objects that we are
interested in, one can obtain $q$-polynomial identities by
counting them in two different ways, and identifying the results. 
One way of counting is by using \lq inclusion-exclusion\rq\
\cite{andrews-red-book}. This was done in \cite{burge}.
We give a detailed proof in Appendix A.

\begin{theorem}\label{GenFunThrm} {\rm \cite{burge}}
\begin{eqnarray*}
\hbox to 5mm{$R(N_1,M_1,N_2,M_2,a,b,\alpha,\beta)$\hss}\\[1mm]
&=& \sum_{k=-\infty}^\infty
q^{k^2(a+b)(\alpha+\beta)+k(\alpha+\beta)(N_1-N_2)+k(a\beta-b\alpha)}\\[-4mm]
&&\qquad\qquad\qquad\qquad\times\quad
g(N_1,M_1,N_2,M_2,k(a+b),k(\alpha+\beta))\\[1mm]
&&\quad-\quad \sum_{k=-\infty}^\infty
q^{k^2(a+b)(\alpha+\beta)+(k(\alpha+\beta)+\alpha)(N_1-N_2)
+k(a\beta+b\alpha+2a\alpha)+a\alpha}\\[-4mm]
&&\qquad\qquad\qquad\qquad\times\quad
g(N_1,M_1,N_2,M_2,k(a+b)+a,k(\alpha+\beta)+\alpha).
\end{eqnarray*}
\end{theorem}

Notice the dependence of the generating function $R$ on {\it two} 
Gaussian polynomials. This is a consequence of the fact that $R$ 
$q$-counts partitions pairs. Also notice its dependence on {\it 
four} finite size parameters. 
Also notice that (taking $N_1=N_2$ and $M_1=M_2$ as we often
will), using (\ref{RochaEq}),
\begin{equation}\label{BRoachEq}
\lim_{\scriptstyle N\to\infty\atop\scriptstyle M\to\infty}
R(N,M,N,M,a,b,\alpha,\beta)=
{1\over (q)_\infty}
\chi^{a+b,\alpha+\beta}_{a,\alpha},
\end{equation}
and thus $R(N_1,M_1,N_2,M_2,a,b,\alpha,\beta)$
is a four-parameter finitised Virasoro character
(up to a factor of $(q)_\infty$).

\subsection{Comparison of notational conventions}

At this point, we wish to compare the different notations used to 
label the Virasoro characters, the corresponding one-point functions, 
and their finite size analogues. It will be sufficient for our purposes 
to compare the notations used in the following four.

\begin{description}
\item[BPZ] 
      In \cite{bpz}, and in the rest of the conformal field theory 
      literature, a Virasoro character is denoted by $\chi^{p, p'}_{r, s}$.
      This notation is used in (\ref{RochaEq}).

\item[FB] 
      In \cite{fb}, a polynomial analogue of  $\chi^{p, p'}_{r, s}$
      that depends on one finite size parameter was obtained (in the 
      process of computing one-point functions). 
      The parameters of Forrester and Baxter are related to those 
      of BPZ as follows:
      $$
      \{p, p', r, s\} \equiv \{K - \mu, K, b, a\}.
      $$
\item[ABBBFV] 
      In \cite{abbbfv}, a polynomial analogue of
      $\chi^{p, p'}_{r, s}$ that depends on one finite size parameter 
      $L$ was obtained. The parameters of ABBBFV are related 
      to those of BPZ as follows:
      $$
      \{p, p', r, s\} \equiv \{\alpha + \beta, K, \beta, \iota\}.
      $$
\item[Burge] 
      In \cite{burge}, a polynomial analogue of $\chi^{p, p'}_{r, s}$ 
      (up to a factor) that depends on two finite size parameters was 
      obtained. The parameters of Burge, and also of this paper, 
      are related to those of BPZ as follows:
      $$
      \{p, p', r, s\} \equiv \{a + b, \alpha + \beta, a, \alpha\}.
      $$
\end{description}

\subsection{The symmetries of the bosonic generating function}

Because we are counting pairs of partitions, rather than single 
ones, the generating function of the partition pairs enjoys a
number of symmetries that are not present in the case of single
partitions. These symmetries are what allows us to use the Burge
transform to obtain more results than in the case of the Bailey
transform. 

The symmetries of the generating function can be expressed in
terms of the following identities between the 
$R(N_1,M_1,N_2,M_2,a,b,\alpha,\beta)$, and 
may be proved directly from Theorem \ref{GenFunThrm}.

\begin{corollary}\label{GenFunCor}
Let $\Delta N=N_1-N_2$, $\Delta M=M_1-M_2$
and $A=N_1M_1+N_2M_2$. Then:
\begin{eqnarray*}
\hbox to 5mm{1.$\quad R(N_1,M_1,N_2,M_2,a,b,\alpha,\beta)$\hss}\\[0.5mm]
&=&R(N_2,M_2,N_1,M_1,b,a,\beta,\alpha);\\[0.5mm]
\hbox to 5mm{2. $\quad R(N_1,M_1,N_2,M_2,a,b,\alpha,\beta)$\hss}\\[0.5mm]
&=&R(M_1,N_1,M_2,N_2,\alpha-\Delta M,\beta+\Delta M,
a+\Delta N,b-\Delta N);\\[0.5mm]
\hbox to 5mm{3.$\quad R(N_1,M_1,N_2,M_2,a,b,\alpha,\beta;q)$\hss}\\[0.5mm]
&=&q^{A}R(N_1,M_1,N_2,M_2,b-\Delta N,a+\Delta N,
\beta+\Delta M,\alpha-\Delta M;q^{-1});\\[0.5mm]
\hbox to 5mm{4.$\quad R(N_1,M_1,N_2,M_2,a,b,\alpha,\beta;q)$\hss}\\[0.5mm]
&=&q^{A}R(N_2,M_2,N_1,M_1,a+\Delta N,b-\Delta N,
\alpha-\Delta M,\beta+\Delta M;q^{-1}).
\end{eqnarray*}
\end{corollary}

\noindent
In fact, these identities may alternatively be deduced combinatorially 
by considering the restricted partition pairs themselves, suitably 
transforming them, and identifying the restriction on the transformed 
pairs. 

The first identity of Corollary \ref{GenFunCor} follows by
interchanging the two partitions that make up the pair.
The second follows by taking the conjugate of each partition.
The third arises from taking the complements of the two partitions:
the first inside a box of dimensions $N_1\times M_1$ and the
second inside a box of dimensions $N_2\times M_2$.
The fourth identity combines the first and third.

\section{The Burge transform}

Let us suppose that we are able to compute the generating function 
of a set of partition pairs that obey certain restrictions in two
different ways. Equating the results, we obtain a $q$-polynomial
identity. Following Burge \cite{burge}, we can use that identity to 
obtain a tree of polynomial identities using a summation formula
that we refer to as {\it The Burge transform}. 

In this section, we recall Burge's derivation of his summation
formula in 3 steps. 

\subsection{A $q$-polynomial form of the $q$-Pfaff-Saalsch\"utz
            summation formula}

The first step is to notice that the $q$-analogue of the 
Pfaff-Saalsch\"utz summation formula (see eq. (3.3.11) of 
\cite{andrews-red-book}, or eq. (4.2) of \cite{gould})
can be written as a $q$-polynomial identity as follows:

\bea\label{baa}
\hbox to 15mm{$\displaystyle
\left[{m_1+n+B\atop n+A}\right]\left[{m_2+n+A\atop n+B}\right]
$\hss}\nonumber\\[0.5mm]
&=& \sum_{i=0}^n q^{(i+A)(i+B)}
\left[{m_1+m_2+n-i\atop n-i}\right]\left[{m_1\atop i+A}\right]
\left[{m_2\atop i+B}\right].
\eea

\noindent
The second step is to use the above identity to prove the
following lemma 

\begin{lemma}\label{BurgeSumLem}
Let $\Delta N=N_1-N_2$. Then
\begin{eqnarray*}
\hbox to 5mm{$\displaystyle
g(N_1,M_1,N_2,M_2,x,y)
$\hss}\\[0.5mm]
&=&\sum_{i=0}^{N_2}
\left( q^{i^2+i\Delta N-x(x+\Delta N)}
\left[{M_1+M_2+N_2-i\atop N_2-i}\right] \right.\\[0.5mm]
&&
\left.\qquad\qquad\phantom{M\atop N}\times
g(i+\Delta N,M_1-i,i,M_2-\Delta N-i,x,y-x)\right).
\end{eqnarray*}
\end{lemma}

\Proof On setting $n=N_2$, $A=\Delta N+x$, $B=-x$,
$m_1=M_1+\Delta N+2x-y$, and $m_2=M_2-\Delta N-2x+y$ in (\ref{baa}),
we obtain:
\begin{eqnarray*}
\hbox to 10mm{$\displaystyle
\left[{N_1+M_1+x-y\atop N_1+x}\right]
\left[{N_2+M_2-x+y\atop N_2-x}\right]
$\hss}\\[0.5mm]
&=&
\sum_{i=0}^{N_2}
\left( q^{(i-x)(i+x+\Delta N)}
\left[{M_1+M_2+N_2-i\atop N_2-i}\right] \right.\\[0.5mm]
&&\qquad\qquad\times
\left.
\left[{M_1+\Delta N+2x-y\atop \Delta N+x+i}\right]
\left[{M_2-\Delta N-2x+y\atop i-x}\right] \right),
\end{eqnarray*}
which is the desired result.
\cqfd

\subsection{The Burge summation formula}

Finally, we are in a position to prove the Burge transform, or
equivalently, the Burge summation formula.

\begin{theorem}\label{BurgeSumThrm}
Let $\Delta N=N_1-N_2$, $\Delta M=M_1-M_2$,
$-a\le\Delta N\le b$, $-\beta\le\Delta M\le\alpha$ and
$b-\beta\le\Delta M+\Delta N\le \alpha-a$.
Then
\begin{eqnarray*}
\hbox to 5mm{$\displaystyle
R(N_1,M_1,N_2,M_2,a,b,\alpha,\beta)
$\hss}\\[0.5mm]
&=&\sum_{i=0}^{N_2}
\left( q^{i^2+i\Delta N}
\left[{M_1+M_2+N_2-i\atop N_2-i}\right]\right.\\[0.5mm]
&&
\left.\qquad\phantom{M\atop N}\times
R(i+\Delta N,M_1-i,i,M_2-\Delta N-i,a,b,\alpha-a,\beta-b)
\right).
\end{eqnarray*}
\end{theorem}

\Proof 
Following Burge \cite{burge}, this expression is proved by applying 
Lemma \ref{BurgeSumLem} to each term of the generating function given 
in Theorem \ref{GenFunThrm}. For the $k$th term of the first summation 
of the expression for $R(N_1,M_1,N_2,M_2,a,b,\alpha,\beta)$, we take 
$x=k(a+b)$ and $y=k(\alpha+\beta)$ in Lemma \ref{BurgeSumLem} to obtain:
\begin{eqnarray*}
\hbox to 5mm{$\displaystyle
q^{k(\alpha+\beta)(k(a+b)+\Delta N)+k(a\beta-b\alpha)}
g(N_1,M_1,N_2,M_2,k(a+b),k(\alpha+\beta))
$\hss}\\[0.5mm]
&=&
q^{k(\alpha+\beta)(k(a+b)+\Delta N)+k(a\beta-b\alpha)}\\[0.5mm]
&&
\times
\sum_{i=0}^{N_2}
\left(
q^{i^2+i\Delta N-k(a+b)(k(a+b)+\Delta N)}
\left[{M_1+M_2+N_2-i\atop N_2-i}\right] \right.\\[0.5mm]
&&
\left.\phantom{M\atop N}\!\times
g(i+\Delta N,M_1-i,i,M_2-\Delta N-i,k(a+b),k(\alpha+\beta-a-b))
\right)\\[0.5mm]
&=&
\sum_{i=0}^{N_2}
\left(
q^{i^2+i\Delta N}
\left[{M_1+M_2+N_2-i\atop N_2-i}\right]
q^{k(\alpha+\beta-a-b)(k(a+b)+\Delta N)+k(a\beta-b\alpha)}
\right.\\[0.5mm]
&&
\left.\phantom{M\atop N}\!\times
g(i+\Delta N,M_1-i,i,M_2-\Delta N-i,k(a+b),k(\alpha+\beta-a-b))
\right)\!.
\end{eqnarray*}
Similarly, for the $k$th term of the second summation of the
expression for $R(N_1,M_1,N_2,M_2,a,b,\alpha,\beta)$, we
take $x=k(a+b)+a$ and $y=k(\alpha+\beta)+\alpha$ in
Lemma \ref{BurgeSumLem} to obtain:
\begin{eqnarray*}
\hbox to 0mm{$\displaystyle
q^{k(\alpha+\beta)(k(a+b)+\Delta N)+k(a\beta+b\alpha+2a\alpha)
+\alpha(a+\Delta N)}
$\hss}\\[0.0mm]
&&
\qquad\qquad\qquad\qquad\times
g(N_1,M_1,N_2,M_2,k(a+b)+a,k(\alpha+\beta)+\alpha)
\\[2.5mm]
&=&
q^{k(\alpha+\beta)(k(a+b)+\Delta N)+k(a\beta+b\alpha+2a\alpha)
+\alpha(a+\Delta N)}\\[0.5mm]
&&
\times\sum_{i=0}^{N_2}
\left(
q^{i^2+i\Delta N-(k(a+b)+a)(k(a+b)+a+\Delta N)}
\left[{M_1+M_2+N_2-i\atop N_2-i}\right]
\right.\\[0.5mm]
&&
\left.\phantom{M\atop N}\times
g(i+\Delta N,M_1-i,i,M_2-\Delta N-i,
\right.\\[-2mm]
&&
\left.\qquad\qquad\qquad\qquad\qquad\qquad\phantom{M\atop N}
k(a+b)+a,k(\alpha+\beta-a-b)+\alpha-a)
\right)\\[0.5mm]
&=&
\sum_{i=0}^{N_2}
\left(
q^{i^2+i\Delta N}
\left[{M_1+M_2+N_2-i\atop N_2-i}\right] \right.\\[0.5mm]
&&
\phantom{M\atop N}\times
q^{k(\alpha+\beta-a-b)(k(a+b)+\Delta N)
+k(a(\beta-b)+b(\alpha-a)+2a(\alpha-a))+(\alpha-a)(a+\Delta N)}
\\[0.5mm]
&&
\left.\phantom{M\atop N}\times
g(i+\Delta N+i,M_1-i,i,M_2-\Delta N-i,
\right.\\[-2mm]
&&
\left.\qquad\qquad\qquad\qquad\qquad\qquad\phantom{M\atop N}
k(a+b)+a,k(\alpha+\beta-a-b)+\alpha-a)
\right),
\end{eqnarray*}
after re-expressing $(k(a+b)+a)(k(a+b)+a+\Delta N)=
k(a+b)(k(a+b)+\Delta N)+k(ab+ba+2a^2)+a(a+\Delta N)$.
Summing over all $k$ for each of these two results, and taking the
difference between the two sums, proves the theorem.
\cqfd

\section{A Burge tree of polynomial identities}

In this section, following Burge \cite{burge}, we obtain an explicit 
algorithm to generate a tree of polynomial identities that depend on 
two finite size parameters. We specialise to the case where $N_1=N_2$ 
and $M_1=M_2$ and so define 
$R(N,M,a,b,\alpha,\beta)=R(N,M,N,M,a,b,\alpha,\beta)$.
We generate a tree-like structure giving fermionic expressions for 
certain $R(N,M,a,b,\alpha,\beta)$, by expressing these certain
$R(N,M,a,b,\alpha,\beta)$ in terms of various 
$R(N^\prime,M^\prime,0,1,1,1)$. The latter form the \lq root\rq\ of the tree, 
and have a particularly simple form which is best obtained via Theorem 
\ref{GenFunThrm}.

\subsection{The root of the tree}

\begin{lemma}\label{RootLem}
$$
\sum_{k=-\infty}^\infty \! q^{k^2+k}
\left(
\left[{N+M\atop N+k}\right] \hskip-3pt \left[{N+M\atop N-k}\right]
-\left[{N+M-1\atop N+k}\right] \hskip-3pt \left[{N+M+1\atop N-k}\right]
\right)
=\delta_{N,0}\,\delta_{M,0}.
$$
\end{lemma}

\Proof In the case where $P>0$, one readily obtains
$$
\left[{P\atop N}\right]
=\left[{P-1\atop N}\right]
+q^{P-N}\left[{P-1\atop N-1}\right],
$$
from (\ref{GaussDef}).
Use of this result in the left side of the premise when
$M+N>0$, yields:
\begin{eqnarray*}
\hbox to 5mm{$\displaystyle
\sum_{k=-\infty}^\infty q^{k^2+k}
\left(
\left[{N+M\atop N-k}\right] 
\left(
\left[{N+M-1\atop N+k}\right]+q^{M-k} \left[{N+M-1\atop N+k-1}\right]
\right)\right.
$\hss}\\[0.5mm]
&&\qquad\qquad\quad\left.
-\left[{N+M-1\atop N+k}\right]
\left( \left[{N+M\atop N-k}\right]+q^{M+k+1} \left[{N+M\atop N-k-1}\right]
\right)
\right)
\\[0.5mm]
&=&
\sum_{k=-\infty}^\infty q^{k^2+M}
\left[{N+M\atop N-k}\right] \left[{N+M-1\atop N+k-1}\right]\\[0.5mm]
&&\qquad\qquad\qquad
-\sum_{k=-\infty}^\infty q^{k^2+2k+1+M}
\left[{N+M\atop N-k-1}\right] \left[{N+M-1\atop N+k}\right]
\\[0.5mm]
&=&
\sum_{k=-\infty}^\infty q^{k^2+M}
\left[{N+M\atop N-k}\right] \left[{N+M-1\atop N+k-1}\right]\\[0.5mm]
&&\qquad\qquad\qquad
-\sum_{k=-\infty}^\infty q^{k^2+M}
\left[{N+M\atop N-k}\right] \left[{N+M-1\atop N+k-1}\right],
\end{eqnarray*}
having shifted the second summation index $k\rightarrow k-1$.
The result is manifestly $0$.
When $M=N=0$, only the $k=0$ term of the first summation is non-zero.
Its value is $1$ which proves the lemma.
\cqfd

\begin{corollary}\label{RootCor}
\begin{eqnarray*}
&&1.\ R(N,M,0,1,1,0)\,=\,\delta_{N,0}\,\delta_{M,0};\\[0.5mm]
&&2.\ R(N,M,0,1,1,1)\,=\,\delta_{M,0}.
\end{eqnarray*}
\end{corollary}

\Proof
Substituting $N_1=N_2=N$, $M_1=M_2=M$, $a=\beta=0$ and $b=\alpha=1$
into Theorem \ref{GenFunThrm} results in the left side
of Lemma \ref{RootLem} after changing the sign of the first summation
index. The first result follows.
The second result then also follows on direct application of
Theorem \ref{BurgeSumThrm} with $N_1=N_2=N$, $M_1=M_2=M$,
$a=0$, $b=\alpha=\beta=1$.
\cqfd

The following similar result will not be required
until Section \ref{GenSection}.

\begin{lemma}\label{RootExcpt}
If $b>0$ and $\beta\ge0$ then
\begin{eqnarray*}
R(0,M,0,b,1,\beta)\,=\,\delta_{M,0}.
\end{eqnarray*}
\end{lemma}

\Proof By Theorem \ref{GenFunThrm},
\begin{eqnarray*}
\hbox to 3mm{$R(0,M,0,b,1,\beta)$\hss}\\[1mm]
&=&\!\sum_{k=-\infty}^\infty
q^{k^2b(\beta+1)-kb}
\left[{M+k(b-\beta-1)\atop kb}\right]
\left[{M-k(b-\beta-1)\atop -kb}\right]
\\[5mm]
&&-\sum_{k=-\infty}^\infty
q^{k^2b(\beta+1)+kb}
\left[{M+k(b-\beta-1)-1\atop kb}\right]
\left[{M-k(b-\beta-1)+1\atop -kb}\right]\!.
\end{eqnarray*}
{}From the definition (\ref{GaussDef}) of the Gaussian, when $k\ne0$
each term here contains a zero factor.
Thus
$$
R(0,M,0,b,1,\beta)\,=\,
\left[{M\atop 0}\right] \left[{M\atop 0}\right]
-\left[{M-1\atop 0}\right] \left[{M+1\atop 0}\right].
$$
The lemma then follows since the first term is always 1, whereas
the second term is 0 when $M=0$ and is 1 otherwise.
\cqfd

\subsection{Continued fractions}

A binary tree of fermionic expressions for certain
$R(N,M,a,b,\alpha,\beta)$ may now be obtained via
Corollary \ref{RootCor}
using the notion of a continued fraction.

Let $s$ and $r$ be positive coprime integers, or $s=1$ and $r=0$.
The {\it continued fraction} $(c_1,c_2,c_3,\ldots,c_t)$ for the 
pair $\cf sr$ is defined as follows. The continued fraction for 
$\cf10$ is defined to be the sequence $()$ of length zero.
Then, for other $\cf sr$, the continued fraction
$(c_1,c_2,c_3,\ldots,c_t)$ is defined recursively
by setting $c_1$ to be the largest integer such that
$c_1r\le s$, and taking $(c_2,c_3,c_4,\ldots,c_t)$ to be
the continued fraction for $\cf r{s-c_1r}$.
This recursive procedure clearly terminates.
Note that if $s$ and $r$ are coprime then so are $r$ and $s-c_1r$
(unless $r\le1$).
Also note that if $t>0$ then $c_1=0$ if and only if $r<s$.
Note further that, apart from the cases $\cf sr=\cf11$
and $\cf sr=\cf10$, we obtain $c_t\ge2$.
It is useful to permit the continued fraction for
which $t>1$ and $c_t=1$ and then to equate $(c_1,\ldots,c_{t-1},1)$
and $(c_1,\ldots,c_{t-1}+1)$.

This definition differs from the usual definition of a continued
fraction, but has the advantage of dealing with the
useful additional case $\cf10$.
The connection with the usual definition is given in the following lemma.

\begin{lemma}\label{CFLemma}
Let $s$ and $r$ be positive coprime integers
and $(c_1,c_2,c_3,\ldots,c_t)$ the {\it continued fraction} for $\cf sr$.
Then
$$
{s\over r}=
{c_1+{\displaystyle\strut 1 \over \displaystyle c_2 +
 {\displaystyle\strut 1 \over \displaystyle c_3 +
 {\displaystyle\strut 1 \over {\lower-5pt\hbox{$\vdots$}
 \over \displaystyle\strut c_{t-1} +
 {\displaystyle\strut 1 \over \displaystyle c_t}}}}}}
$$
\end{lemma}

\Proof
Clearly $t\ge1$ and the result holds in the case $t=1$.
For the purposes of induction, assume that the result
holds for arbitrary $t\ge1$.
For $t\ge1$, let $\cf sr$ be such as to have continued fraction
$(c_1,c_2,c_3,\ldots,c_{t+1})$.
Then, by definition, $(c_2,c_3,\ldots,c_{t+1})$ is the continued
fraction for $\cf r{s-c_1r}$.
The induction hypothesis now implies that the quotient
$r/(s-c_1r)$ is given by the right side
of the premise once each subscript is increased by one.
The result at $t+1$ now follows because $s/r=c_1+1/(r/(s-c_1r))$.
The lemma is then proved by induction.
\cqfd

\subsection{A \lq Bailey-type\rq\ algorithm to generate a Burge tree}

\begin{lemma}\label{OnePassLem}
For $c\ge0$,
\begin{eqnarray*}
\hbox to 5mm{$\displaystyle
R(N,M,a,b,\alpha+ca,\beta+cb)
$\hss}\\[0.5mm]
&=&\hskip-5mm
\sum_{\scriptstyle n\:{\rm in}\:(c,N)\atop\scriptstyle\wt(n)\,\le\,M}
\left(
q^{n_1^2+n_2^2+\cdots+n_c^2}
\prod_{i=1}^{c}
\left[{2\left(M-\sum_{j=1}^i n_j\right)+n_{i-1}+n_i\atop n_{i-1}-n_i}\right]
\right.\\[0.5mm]
&&\left.\hskip45mm\times\phantom{M\atop N}
R(n_c,M-\wt(n),a,b,\alpha,\beta)\right),
\end{eqnarray*}
\end{lemma}
the sum being over all partitions $n=(n_1,n_2,\ldots,n_c)$ in $(c,N)$
for which $\wt(n)\le M$.

\Proof The result holds trivially when $c=0$.
For the purposes of induction, assume that the result holds
for a fixed $c\ge0$.
Using this followed by Theorem \ref{BurgeSumThrm} (with $N_1=N_2=n_c$
and $M_1=M_2=M-\wt(n)$) gives:
\begin{eqnarray*}
\hbox to 5mm{$\displaystyle
R(N,M,a,b,\alpha+(c+1)a,\beta+(c+1)b)
$\hss}\\[0.5mm]
&=&\hskip-5mm
\sum_{\scriptstyle n\:{\rm in}\:(c,N)\atop\scriptstyle\wt(n)\,\le\, M}
\left(
q^{n_1^2+n_2^2+\cdots+n_c^2}
\prod_{i=1}^{c}
\left[{2\left(M-\sum_{j=1}^i n_j\right)+n_{i-1}+n_i\atop n_{i-1}-n_i}\right]
\right.\\[0.5mm]
&&\left.\hskip25mm\times\phantom{M\atop N}
R(n_c,M-\wt(n),a,b,\alpha+a,\beta+b)\right)\\[0.5mm]
&=&\hskip-5mm
\sum_{\scriptstyle n\:{\rm in}\:(c,N)\atop\scriptstyle\wt(n)\,\le\, M}
\left(
q^{n_1^2+n_2^2+\cdots+n_c^2}
\prod_{i=1}^{c}
\left[{2\left(M-\sum_{j=1}^i n_j\right)+n_{i-1}+n_i\atop n_{i-1}-n_i}\right]
\right.\\[0.5mm]
&&\left.\hskip25mm\times\phantom{M\atop N}
\sum_{n_{c+1}=0}^{n_c} q^{n_{c+1}}
\left[{2\left(M-\sum_{j=1}^c n_j\right)+n_{c}-n_{c+1}\atop n_{c}-n_{c+1}}\right]
\right.\\[0.5mm]
&&\left.\hskip25mm\times\phantom{M\atop N}
R(n_{c+1},M-\wt(n)-n_{c+1},a,b,\alpha+a,\beta+b)\right)\\[0.5mm]
&=&\hskip-5mm
\sum_{\scriptstyle n\:{\rm in}\:(c+1,N)\atop\scriptstyle\wt(n)\,\le\, M}
\left(
q^{n_1^2+n_2^2+\cdots+n_c^2+n_{c+1}^2}
\prod_{i=1}^{c+1}
\left[{2\left(M-\sum_{j=1}^i n_j\right)+n_{i-1}+n_i\atop n_{i-1}-n_i}\right]
\right.\\[0.5mm]
&&\left.\hskip45mm\times\phantom{M\atop N}
R(n_{c+1},M-\wt(n),a,b,\alpha,\beta)\right).
\end{eqnarray*}
Thus the result holds at $c+1$, whereupon the theorem follows by induction.
\cqfd

\subsection{A \lq Burge-type\rq\ algorithm to generate a Burge tree}

\begin{definition}
The following and subsequent theorems give fermionic expressions
involving a summation taken over various sequences of partitions
that satisfy certain constraints.
Typically, given a sequence of integers $(c_1,c_2,\ldots,c_t)$,
the sum will be over all sequences  of partitions
$n^{(1)}$, $n^{(2)},\ldots,n^{(x)}$, for $x\le t$
(where $x$ may or may not be fixed) where each partition $n^{(k)}$
is in $(c_k,n_{k,0})$ for a certain integer $n_{k,0}$ and
$\wt(n^{(k)})\le w_k$ for a certain integer $w_k$.
The values of $n_{k,0}$ and $w_k$ will be specified recursively.
The parts of $n^{(k)}$ will be denoted $n_{k,i}$ for $1\le k\le c_k$
and listed (conventionally) in non-increasing order. Thus:
$$
n_{k,0}\ge n_{k,1}\ge n_{k,2}\ge\cdots\ge n_{k,c_k}\ge 0
$$
for $1\le k\le x$.

In later results, some of the $n_{k,0}$ or the $w_k$ will be
specified to be $\infty$. This will naturally mean that the
largest part of $n^{(k)}$ is unbounded or that $\wt(n^{(k)})$ is
unbounded respectively. Also note that either of these values
might be specified $\infty\pm a$ for some integer $a$.
This sum/difference should also be taken to be $\infty$.
\cqfd
\end{definition}

The following rather technical result will be used as a stepping
stone to later theorems.

\begin{lemma}\label{CFThrm}
Let $p$ and $p^\prime$ be positive coprime integers and let
$(c_1,c_2,\ldots,c_t)$ be the continued fraction for $\cf{p^\prime}p$
satisfying $t\ge1$.
For $1\le u\le t$ and $v\ge u-1$,
let $r$ and $s$ be such that
$(c_1,c_2,\ldots,c_{u-1},d_{u},d_{u+1},\ldots,d_v)$ is the
continued fraction for $\cf sr$.
Now let $\hat p$ and $\hat p^\prime$ be such that
$(c_{u},c_{u+1},\ldots,c_t)$ is the continued fraction
for $\cf{\hat p^\prime}{\hat p}$, and
let $\hat r$ and $\hat s$ be such that
$(d_{u},d_{u+1},\ldots,d_v)$ is the continued fraction
for $\cf{\hat s}{\hat r}$.
Then $R(N,M,r,p-r,s,p^\prime-s)$ may be expressed in terms of
$R(N^\prime,M^\prime,\hat r,\hat p-\hat r,\hat s,\hat p^\prime-\hat s)$,
for various $N^\prime$ and $M^\prime$,
as follows:
\begin{eqnarray*}
\hbox to 0mm{$\displaystyle
R(N,M,r,p-r,s,p^\prime-s)
$\hss}\\[0.5mm]
\hskip-3mm
&=&
\sum
\left(
q^{\sum_{k=1}^{u-1}\sum_{i=1}^{c_k} n_{k,i}^2}
\prod_{k=1}^{u-1}\prod_{i=1}^{c_k}\!
\left[{2\left(w_k-\sum_{j=1}^i n_{k,j}\right)\!+n_{k,i-1}+n_{k,i}
\atop n_{k,i-1}-n_{k,i}}\right]\!\right.
\\[0.5mm]
&&\left.\hskip45mm\times\phantom{M\atop N}
R(n_{u,0}, w_u, \hat r,\hat p-\hat r,\hat s,\hat p^\prime-\hat s)
\right),
%
\end{eqnarray*}
where the sum is over all sequences
$n^{(1)}$, $n^{(2)},\ldots,n^{(u-1)}$ of partitions for which, for $1\le k<u$,
the partition $n^{(k)}$ is in $(c_k,n_{k,0})$ and satisfies
$\wt(n^{(k)})\le w_k$, where we define $w_1=M$, $w_k=n_{k-1,c_{k-1}}$
for $2\le k\le u$, $n_{1,0}=N$, and $n_{k,0}=w_{k-1}-\wt(n^{(k-1)})$
for $2\le k\le u$.
\end{lemma}

\Proof If $u=1$, $t>0$ and $v\ge0$ then the expression holds trivially.
We shall proceed by induction on $u$ keeping the differences
$t-u$ and $v-u$ fixed.

For the purposes of induction, assume that the result holds for a
fixed $u\ge1$.
Now let $\cf{p^\prime}p$ have continued fraction $(c_1,c_2,\ldots,c_{t+1})$,
$\cf sr$ have continued fraction $(c_1,c_2,\ldots,c_{u})$,
$\cf{\hat p^\prime}{\hat p}$ have continued fraction
$(c_{u+1},c_{u+2},\ldots,c_{t+1})$,
and $\cf{\hat s}{\hat r}$ have continued fraction
$(d_{u+1},d_{u+2},\ldots,d_{v+1})$,
Then, by definition,
$\cf p{p^\prime-c_1p}$ has continued fraction $(c_2,c_3,\ldots,c_{t+1})$
and $\cf r{s-c_1r}$ has continued fraction $(c_2,c_3,\ldots,c_{u})$.
Thereupon, on using Lemma \ref{OnePassLem} followed by Corollary
\ref{GenFunCor}(2) and then the induction hypothesis,
\begin{eqnarray*}
\hbox to 5mm{$\displaystyle
R(N,M,r,p-r,s,p^\prime-s)
$\hss}\\[0.5mm]
&=&\hskip-4mm
\sum_{\lower5mm\hbox{$\scriptstyle
n^{(1)}\:{\rm in}\:(c_1,n_{1,0})\atop\scriptstyle\wt(n^{(1)})\,\le\,w_1$}}
\hskip-4mm
\left(q^{n_{1,1}^2+n_{1,2}^2+\cdots+n_{1,c_1}^2}
\prod_{i=1}^{c_1}
\left[{2\left(w_1-\sum_{j=1}^i n_j\right)+n_{i-1}+n_i\atop n_{i-1}-n_i}\right]
\right.\\[-2mm]
&&\left.\hskip22mm\phantom{M\atop N}\times
R(w_2,n_{2,0},r,p-r,s-c_1r,p^\prime-s-c_1p+c_1r)
\right)\\[0.5mm]
&=&\hskip-4mm
\sum_{\lower5mm\hbox{$\scriptstyle
n^{(1)}\:{\rm in}\:(c_1,n_{1,0})\atop\scriptstyle\wt(n^{(1)})\,\le\,w_1$}}
\hskip-4mm
\left(q^{n_{1,1}^2+n_{1,2}^2+\cdots+n_{1,c_1}^2}
\prod_{i=1}^{c_1}
\left[{2\left(w_1-\sum_{j=1}^i n_j\right)+n_{i-1}+n_i\atop n_{i-1}-n_i}\right]
\right.\\[-2mm]
&&\left.\hskip22mm\phantom{M\atop N}\times
R(n_{2,0},w_2,s-c_1r,p^\prime-s-c_1p+c_1r,r,p-r)
\right)\\[0.5mm]
&=&\hskip-4mm
\sum_{\lower5mm\hbox{$\scriptstyle
n^{(1)}\:{\rm in}\:(c_1,n_{1,0})\atop\scriptstyle\wt(n^{(1)})\,\le\,w_1$}}
\hskip-4mm
\left(q^{n_{1,1}^2+n_{1,2}^2+\cdots+n_{1,c_1}^2}
\prod_{i=1}^{c_1}
\left[{2\left(w_1-\sum_{j=1}^i n_j\right)+n_{i-1}+n_i\atop n_{i-1}-n_i}\right]
\right)\\[0.5mm]
&&\quad\times\hskip-11mm\sum_{\vbox{\openup1\jot\halign{\hfil
$\scriptstyle #$&$\scriptstyle #$\hfil\cr
\noalign{\vskip1mm}
n^{(2)}\:{\rm in}\:(c_2,n_{2,0}),\,&\wt(n^{(2)})\,\le\,w_2\cr
n^{(3)}\:{\rm in}\:(c_3,n_{3,0}),\,&\wt(n^{(3)})\,\le\,w_3\cr
&\:\vdots\cr
n^{(u)}\:{\rm in}\:(c_u,n_{u,0}),\,&\wt(n^{(u)})\,\le\,w_u\cr}}}
\hskip-14mm
\left(
q^{\sum_{k=2}^{u}\sum_{i=1}^{c_k} n_{k,i}^2}
\phantom{M\atop N}\right.\\[0.5mm]
&&\qquad\qquad\phantom{M\atop N}\times
\prod_{k=2}^{u}\prod_{i=1}^{c_k}
\left[{2\left(w_k-\sum_{j=1}^i n_{k,j}\right)+n_{k,i-1}+n_{k,i}
\atop n_{k,i-1}-n_{k,i}}\right]
\\[0.5mm]
&&\left.\qquad\qquad\phantom{M\atop N}\times
R(n_{u+1,0},w_{u+1},\hat r,\hat p-\hat r,\hat s,\hat p^\prime-\hat s)
\right)\\[0.5mm]
&=&\hskip-11mm\sum_{\vbox{\openup1\jot\halign{\hfil
$\scriptstyle #$&$\scriptstyle #$\hfil\cr
\noalign{\vskip1mm}
n^{(1)}\:{\rm in}\:(c_1,n_{1,0}),\,&\wt(n^{(1)})\,\le\,w_1\cr
n^{(2)}\:{\rm in}\:(c_2,n_{2,0}),\,&\wt(n^{(2)})\,\le\,w_2\cr
&\:\vdots\cr
n^{(u)}\:{\rm in}\:(c_u,n_{u,0}),\,&\wt(n^{(u)})\,\le\,w_u\cr}}}
\hskip-14mm
\left(
q^{\sum_{k=1}^{u}\sum_{i=1}^{c_k} n_{k,i}^2}
\right.\\[0.5mm]
&&\qquad\qquad\phantom{M\atop N}\times
\prod_{k=1}^{u}\prod_{i=1}^{c_k}
\left[{2\left(w_k-\sum_{j=1}^i n_{k,j}\right)+n_{k,i-1}+n_{k,i}
\atop n_{k,i-1}-n_{k,i}}\right]
\\[0.5mm]
&&\left.\qquad\qquad\phantom{M\atop N}\times
R(n_{u+1,0},w_{u+1},\hat r,\hat p-\hat r,\hat s,\hat p^\prime-\hat s)
\right),
\end{eqnarray*}
where $w_{u+1}=n_{u,c_u}$ and $n_{u+1,0}=w_u-\wt(n^{(u)})$.
This is the desired result at $u+1$, $t+1$ and $v+1$.
Hence the lemma is proved by induction.
\cqfd

\begin{theorem}\label{DetOneThrm}
Let $p$ and $p^\prime$ be positive coprime integers and let
$\cf{p^\prime}p$ have continued fraction
$(c_1,c_2,\ldots,c_{t-1},c_t+2)$.
Then let $r$ and $s$ be such that $\cf sr$ has continued fraction
$(c_1,c_2,\ldots,c_{t-1})$.
Then $R(N,M,r,p-r,s,p^\prime-s)$ may be expressed as follows:
\begin{eqnarray*}
\hbox to 0mm{$\displaystyle
R(N,M,r,p-r,s,p^\prime-s)
$\hss}\\[0.5mm]
&&\!=\sum
q^{\sum_{k=1}^t\sum_{i=1}^{c_k} n_{k,i}^2}
\prod_{k=1}^{t}\prod_{i=1}^{c_k}
\left[{2\left(w_k-\sum_{j=1}^i n_{k,j}\right)+n_{k,i-1}+n_{k,i}
\atop n_{k,i-1}-n_{k,i}}\right]
\!,
\end{eqnarray*}
where the sum is over all sequences
$n^{(1)}$, $n^{(2)},\ldots,n^{(t)}$ of partitions for which, for $1\le k\le t$,
the partition $n^{(k)}$ is in $(c_k,n_{k,0})$ and satisfies
$\wt(n^{(k)})\le w_k$, where we define $w_1=M$, $w_k=n_{k-1,c_{k-1}}$
for $2\le k\le t$, $n_{1,0}=N$, and $n_{k,0}=w_{k-1}-\wt(n^{(k-1)})$
for $2\le k\le t$;
and additionally also satisfy $\wt(n^{(t)})=w_t$.
\end{theorem}

\Proof On using Lemma \ref{CFThrm} with $u=t=v+1$, so that
$\cf{\hat p^\prime}{\hat p}$ has continued fraction $(c_t+2)$ and
$\cf{\hat s}{\hat r}$ has continued fraction $()$, whereupon
$\hat p=1$, $\hat p^\prime=c_t+2$, $\hat r=0$ and $\hat s=1$,
we obtain
\begin{eqnarray*}
\hbox to 5mm{$\displaystyle
R(N,M,r,p-r,s,p^\prime-s)
$\hss}\\[0.5mm]
&=&\hskip-11mm\sum_{\vbox{\openup1\jot\halign{\hfil
$\scriptstyle #$&$\scriptstyle #$\hfil\cr
\noalign{\vskip1mm}
n^{(1)}\:{\rm in}\:(c_1,n_{1,0}),\,&\wt(n^{(1)})\,\le\,w_1\cr
n^{(2)}\:{\rm in}\:(c_2,n_{2,0}),\,&\wt(n^{(2)})\,\le\,w_2\cr
&\:\vdots\cr
n^{(t-1)}\:{\rm in}\:(c_{t-1},n_{t-1,0}),\,&\wt(n^{(t-1)})\,\le\,w_{t-1}\cr}}}
\hskip-20mm
\left(
q^{\sum_{k=1}^{t-1}\sum_{i=1}^{c_k} n_{k,i}^2}
\phantom{M\atop N}\right. \\[0.5mm]
&&\qquad\qquad\times\phantom{M\atop N}
\prod_{k=1}^{t-1}\prod_{i=1}^{c_k}
\left[{2\left(w_k-\sum_{j=1}^i n_{k,j}\right)+n_{k,i-1}+n_{k,i}
\atop n_{k,i-1}-n_{k,i}}\right]
\\[0.5mm]
&&\left.\qquad\qquad\times\phantom{M\atop N}
R(n_{t,0},w_{t},0,1,1,c_t+1)
\right).
\end{eqnarray*}
The use of Lemma \ref{OnePassLem} with $a=0$, $b=\alpha=\beta=1$,
$c=c_t$, $N=n_{t,0}$, $M=w_t$,
followed by an application of Corollary \ref{RootCor}, results in:
\begin{eqnarray*}
\hbox to 5mm{$\displaystyle
R(n_{t,0},w_t,0,1,1,c_t+1)
$\hss}\\[0.5mm]
&=&
\hskip-3mm
\sum_{\lower5mm\hbox{$\scriptstyle
n^{(t)}\:{\rm in}\:(c_t,n_{t,0})\atop\scriptstyle\wt(n^{(t)})\,\le\,w_t$}}
\hskip-5mm
\left(
q^{n_{t,1}^2+n_{t,2}^2+\cdots+n_{t,c_t}^2}
\prod_{i=1}^{c_t}
\left[{2\left(w_t-\sum_{j=1}^i n_{t,j}\right)+n_{t,i-1}+n_{t,i}\atop
  n_{t,i-1}-n_{t,i}}\right]
\right.\\[-2mm]
&&\left.\hskip40mm\times\phantom{M\atop N}
R(n_{t,c_t},w_t-\wt(n^{(t)}),0,1,1,1)\right)\\[0.5mm]
&=&
\hskip-3mm
\sum_{\lower5mm\hbox{$\scriptstyle
n^{(t)}\:{\rm in}\:(c_t,n_{t,0})\atop\scriptstyle\wt(n^{(t)})\,\le\,w_t$}}
\hskip-5mm
\left(
q^{n_{t,1}^2+n_{t,2}^2+\cdots+n_{t,c_t}^2}
\prod_{i=1}^{c_t}
\left[{2\left(w_t-\sum_{j=1}^i n_{t,j}\right)+n_{t,i-1}+n_{t,i}\atop
  n_{t,i-1}-n_{t,i}}\right]
\right.\\[-2mm]
&&\left.\hskip40mm\times\phantom{M\atop N}
\delta_{w_t,\wt(n^{(t)})}
\right).
\end{eqnarray*}
which, when substituted in the previous expression, proves the theorem.
\cqfd

In fact, the particular $r$ and $s$ specified in Theorem \ref{DetOneThrm}
have a simple characterisation. We first prove the following:

\begin{lemma}\label{WhichLem}
Let $p$ and $p^\prime$ be positive coprime integers and let
$(c_1,c_2,\ldots,c_t)$ be the continued fraction for $\cf {p^\prime}p$
satisfying $c_t\ge2$. Let $r_0$ and $s_0$ be such that
$(c_1,c_2,\ldots,c_{t-1})$ is a continued fraction for
$\cf{s_0}{r_0}$.
Then:
\begin{enumerate}
\item
$p^\prime r_0-ps_0=(-1)^t$;
\item
$r_0\le p/2$ and $s_0\le p^\prime/2$.
\item
$(c_1,c_2,\ldots,c_{t-1},c_t-1)$ is a continued fraction
for $\cf{p^\prime-s_0}{p-r_0}$.
\end{enumerate}
\end{lemma}

\Proof
For each part, the $t=1$ case is immediate.
For the purposes of induction,
assume that each part holds for a fixed $t\ge1$.
Now let $\cf {p^\prime}p$ and $\cf{s_0}{r_0}$ have continued fractions
$(c_1,c_2,\ldots,c_{t+1})$ and $(c_1,c_2,\ldots,c_{t})$
respectively.
Then, by definition, $\cf p{p^\prime-c_1p}$ and $\cf{r_0}{s_0-c_1r_0}$
have continued fractions
$(c_2,c_3,\ldots,c_{t+1})$ and $(c_2,c_3,\ldots,c_{t})$
respectively.
Then, by the induction hypothesis for the first part,
$p(s_0-c_1r_0)-(p^\prime-c_1p)r_0=(-1)^t$, whereby
$p^\prime r_0-ps_0=(-1)^{t+1}$.
The first part now follows by induction.

By the induction hypothesis for the second part,
$r_0\le p/2$ and $s_0-c_1r_0\le(p^\prime-c_1p)/2$ whereupon
$s_0\le p^\prime/2+c_1(r_0-p/2)\le p^\prime/2$.
The second part now also follows by induction.

By the induction hypothesis for the third part,
$\cf{p-r_0}{p^\prime-c_1p-s_0+c_1r_0}$ has continued fraction
$(c_2,c_3,\ldots,c_{t+1}-1)$.
Thereupon $(c_1,c_2,c_3,\ldots,c_{t+1}-1)$ is a continued fraction for
$\cf{p^\prime-s_0}{p-r_0}$.
The third part now also follows by induction.
\cqfd

\noindent
It is now a straightforward task to show that the $r_0$ and $s_0$ specified
in this lemma are the smallest positive integers for which
$\vert\, p^\prime r_0-ps_0\,\vert=1$.
Thus these are the values denoted $r_{\rm min}$ and $s_{\rm min}$
in \cite{bm}.

\begin{note}\label{Anote} {\rm
By using the first identity of Corollary \ref{GenFunCor},
the fermionic expression for $R(N,M,p-r,r,p^\prime-s,s)$ is
immediately obtained from that of $R(N,M,r,p-r,s,p^\prime-s)$.
In the case of $r_0$ and $s_0$, we thus also have, by
virtue of Theorem \ref{DetOneThrm}, the fermionic
expression for $R(N,M,r,p-r,s,p^\prime-s)$ where
$r=p-r_0$ and $s=p^\prime-s_0$.
Thus for specific $p$ and $p'$ with $\cf {p'}p$ having
continued fraction $(c_1,c_2,\ldots,c_t)$, we already have
fermionic expressions for
$R(N,M,r,p-r,s,p^\prime-s)$ for two different pairs of
$r$ and $s$: namely those for which $\cf sr$ has continued
fractions $(c_1,c_2,\ldots,c_{t-1})$ and
$(c_1,c_2,\ldots,c_t-1)$.

In fact, it is possible to calculate
$R(N,M,p-r_0,r_0,p^\prime-s_0,s_0)$ directly using the
technique of this section, but starting (instead)
from $R(N^\prime,M^\prime,1,0,1,1)=\delta_{M^\prime,0}$.
This also results in precisely the same expression for
$R(N,M,p-r_0,r_0,p^\prime-s_0,s_0)$ as for
$R(N,M,r_0,p-r_0,s_0,p^\prime-s_0)$.}
\end{note}

\section{An extended Burge tree}\label{GenSection}

In this section, we deduce a further relationship between the sets
of partition pairs. This is then exploited to provide fermionic
expressions for further $R(N,M,r,p-r,s,p^\prime-s)$.
Once again, the continued fraction of $(p^\prime:p)$ determines
for which $r$ and $s$, the fermionic expression for
$R(N,M,r,p-r,s,p^\prime-s)$ may be determined.
However here, for particular $p$ and $p'$, this may be achieved
for various $r$ and $s$.

\begin{lemma}\label{ShiftLemma}
Let $M$, $N$, $b$, and $\beta$ be non-negative integers.
Then
\begin{eqnarray*}
&&1.\ R(N,M+1,1,b,0,\beta+1)\,=\,q^N R(N,M,1,b,1,\beta);\\[0.5mm]
&&2.\ R(N+1,M,0,b+1,1,\beta)\,=\,q^M R(N,M,1,b,1,\beta);\\[0.5mm]
&&3.\ R(N+1,M,0,b+1,1,\beta)\,=\,q^{M-N} R(N,M+1,1,b,0,\beta+1);\\[0.5mm]
&&4.\ R(N+1,M,0,b+1,1,\beta)\,=\,q^{M-N} R(M+1,N,0,\beta+1,1,b).
\end{eqnarray*}
\end{lemma}

\Proof
The first expression is proved by setting up a bijection between the
sets ${\cal R}(N,M+1,N,M+1,1,b,0,\beta+1)$ and
${\cal R}(N,M,N,M,1,b,1,\beta)$.
Let $(q,p)\in{\cal R}(N,M+1,N,M+1,1,b,0,\beta+1)$ whence from
(\ref{restrictions}),
\begin{displaymath}
\mathalign{
p_i\:&\ge q_{i}+1&\qquad1\le i\le N,\cr
q_i\:&\ge p_{i+1-b}-\beta&\qquad b\le i\le N.\cr}
\end{displaymath}
Since each $q_i\ge0$ and each $p_i\le M+1$, the first of these
inequalities implies that $p_i\ge1$ and $q_i\le M$.
Thereupon, on setting $q_i^\prime=q_i$ and $p_i^\prime=p_i-1$
for $1\le i\le N$, we have:
\begin{displaymath}
\mathalign{
p_i^\prime\:&\ge q_{i}^\prime&\qquad1\le i\le N,\cr
q_i^\prime\:&\ge p_{i+1-b}^\prime+1-\beta&\qquad b\le i\le N.\cr}
\end{displaymath}
This ensures that
$(q^\prime,p^\prime)\in{\cal R}(N,M,N,M,1,b,1,\beta)$.
The map $(q,p)\mapsto(q^\prime,p^\prime)$ clearly defines
a bijection between ${\cal R}(N,M+1,N,M+1,1,b,0,\beta+1)$
and ${\cal R}(N,M,N,M,1,b,1,\beta)$.
Thereupon, since $\wt((q^\prime,p^\prime))=\wt((q,p))-N$,
the first expression is proved.
(The first expression may also be proved via the expression
(\ref{GenFunThrm}) for the generating function.)

The second expression follows from the first by transforming both
sides using Corollary \ref{GenFunCor}(2).
Thereupon, the third expression results from combining these two.
A further application of Corollary \ref{GenFunCor}(2) then yields
the fourth expression.
\cqfd

\medskip
As we shall see, this lemma enables us to determine fermionic
expressions for $R(N,M,r,p-r,s,p^\prime-s)$ where the possible
$r$ and $s$ are determined by $p$ and $p^\prime$ as follows.
Let $\cf{p^\prime}p$ have continued fraction
$(c_1,c_2,\ldots,c_t+2)$ with $t\ge1$.
Then we may obtain the fermionic expression for
$R(N,M,r,p-r,s,p^\prime-s)$ if $\cf sr$ has continued fraction
$(c_1,c_2,\ldots,c_{u-1},d)$ where $1\le u\le t$ and $1\le d\le c_{u}$.
(Here we include the possibility that $u>1$ and $d=1$.
In such a case, the continued fractions
$(c_1,c_2,\ldots,c_{u-1},1)$ and $(c_1,c_2,\ldots,c_{u-1}+1)$ are
equated.)
This is done in Theorem \ref{GenThrm}.
Of course, for these $r$ and $s$, Corollary \ref{GenFunCor}(1) then also
yields $R(N,M,p-r,r,p^\prime-s,s)$.

In the following result, we restrict to the $N>0$ case, the $N=0$ case
having been already dealt with in Lemma \ref{RootExcpt}.

\begin{theorem}\label{DetGenThrm}
Let $p$ and $p^\prime$ be positive coprime integers and let
$\cf{p^\prime}p$ have continued fraction
$(c_1,c_2,\ldots,c_{t-1},c_t+2)$ with $t\ge1$.
Then $R(N,M,0,p,1,p^\prime-1)$ may be expressed as follows
when $N>0$:
\begin{eqnarray*}
\hbox to 0mm{$\displaystyle
R(N,M,0,p,1,p^\prime-1)
$\hss}\\[0.5mm]
&=&\!\!\!q^M \sum
q^{\sum_{k=1}^x \sum_{i=1}^{c_k} n_{k,i}(n_{k,i}-1)}
\prod_{k=1}^{x}\prod_{i=1}^{c_k}
\left[{2\left(w_k\!-\!\sum_{j=1}^i n_{k,j}\right)\!+\!n_{k,i-1}\!+\!n_{k,i}
\atop n_{k,i-1}\!-\!n_{k,i}}\right]
\!\!,
\end{eqnarray*}
where the sum is over all sequences
$n^{(1)}$, $n^{(2)},\ldots,n^{(x)}$ of partitions with $1\le x\le t$ for which,
for $1\le k\le x$,
the partition $n^{(k)}$ is in $(c_k,n_{k,0})$ and satisfies
$\wt(n^{(k)})\le w_k$ with $\wt(n^{(x)})=w_x$,
where we define $w_1=M$, $w_k=n_{k-1,c_{k-1}}-1$
for $2\le k\le x$, $n_{1,0}=N$, and $n_{k,0}=w_{k-1}+1-\wt(n^{(k-1)})$
for $2\le k\le x$;
and additionally also
satisfies $n_{k,c_k}>0$ for $1\le k<x$ and
satisfies $n_{x,c_x}=0$ whenever $x<t$.
\end{theorem}

\Proof If $t=1$ then necessarily $p^\prime=c_1+2$ and $p=1$.
Then, with $n_{1,0}=N$ and $w_1=M$, Lemma \ref{OnePassLem} implies that
\begin{eqnarray*}
\hbox to 5mm{$\displaystyle
R(N,M,0,1,1,c_t+1)
$\hss}\\[0.5mm]
&=&
\hskip-3mm
\sum_{\lower5mm\hbox{$\scriptstyle
n^{(1)}\:{\rm in}\:(c_1,n_{1,0})\atop\scriptstyle\wt(n^{(1)})\,\le\,w_1$}}
\hskip-6mm
\left(
q^{\sum_{i=1}^{c_1} n_{1,i}^2}
\prod_{i=1}^{c_1}
\left[{2\left(w_1-\sum_{j=1}^i n_{1,j}\right)+n_{1,i-1}+n_{1,i}\atop
  n_{1,i-1}-n_{1,i}}\right]
\right.\\[-2mm]
&&\left.\hskip40mm\times\phantom{M\atop N}
R(n_{1,c_1},w_1-\wt(n^{(1)}),0,1,1,1)\right)\\[0.5mm]
&=&
\hskip-3mm
\sum_{\lower5mm\hbox{$\scriptstyle
n^{(1)}\:{\rm in}\:(c_1,n_{1,0})\atop\scriptstyle\wt(n^{(1)})\,\le\,w_1$}}
\hskip-6mm
\left(
q^{\sum_{i=1}^{c_1} n_{1,i}^2}
\prod_{i=1}^{c_1}
\left[{2\left(w_1-\sum_{j=1}^i n_{1,j}\right)+n_{1,i-1}+n_{1,i}\atop
  n_{1,i-1}-n_{1,i}}\right]
\right.\\[-2mm]
&&\left.\hskip40mm\times\phantom{M\atop N}
\delta_{w_1,\wt(n^{(1)})}
\right).
\end{eqnarray*}
on using Corollary \ref{RootCor}.
Therefore, we require those partitions $n^{(1)}$ for which $\wt(n^{(1)})=w_1=M$.
Furthermore, we then have that
$\sum_{i=1}^{c_1} n_{1,i}^2 = M+\sum_{i=1}^{c_t} n_{1,i}(n_{1,i}-1)$,
and therefore the theorem holds in the case $t=1$, when $x$ can
only take the value $x=1$.

For the purposes of induction, assume that the result holds for a
fixed $t\ge1$.
Now let $\cf{p^\prime}p$ have continued fraction $(c_1,c_2,\ldots,c_{t+1}+2)$.
Then, by definition,
$\cf p{p^\prime-c_1p}$ has continued fraction $(c_2,c_3,\ldots,c_{t+1}+2)$.
Thereupon, on using Lemma \ref{OnePassLem},
\begin{eqnarray*}
\hbox to 5mm{$\displaystyle
R(N,M,0,p,1,p^\prime-1)
$\hss}\\[0.5mm]
&=&
\hskip-3mm
\sum_{\lower5mm\hbox{$\scriptstyle
n^{(1)}\:{\rm in}\:(c_1,n_{1,0})\atop\scriptstyle\wt(n^{(1)})\,\le\,w_1$}}
\hskip-6mm
\left(
q^{\sum_{i=1}^{c_1} n_{1,i}^2}
\prod_{i=1}^{c_1}
\left[{2\left(w_1-\sum_{j=1}^i n_{1,j}\right)+n_{1,i-1}+n_{1,i}\atop
  n_{1,i-1}-n_{1,i}}\right]
\right.\\[-2mm]
&&\left.\hskip30mm\times\phantom{M\atop N}
R(n_{1,c_1},w_1-\wt(n^{(1)}),0,p,1,p^\prime-c_1p-1)\!\right).
\end{eqnarray*}
For those cases where $n_{1,c_1}=0$, by virtue of Lemma \ref{RootExcpt},
we obtain:
\begin{eqnarray*}
\hbox to 5mm{$\displaystyle
R(N,M,0,p,1,p^\prime-1)
$\hss}\\[0.5mm]
&=&
\hskip-3mm
\sum_{\lower5mm\hbox{$\scriptstyle
n^{(1)}\:{\rm in}\:(c_1,n_{1,0})\atop\scriptstyle\wt(n^{(1)})\,=\,w_1$}}
\hskip-6mm
\left(
q^{\sum_{i=1}^{c_1} n_{1,i}^2}
\prod_{i=1}^{c_1}
\left[{2\left(w_1-\sum_{j=1}^i n_{1,j}\right)+n_{1,i-1}+n_{1,i}\atop
  n_{1,i-1}-n_{1,i}}\right]
\right.\\[-2mm]
&&\left.\hskip30mm\times\phantom{M\atop N}
\delta_{w_1,\wt(n^{(1)})}\right)\\[0.5mm]
&=&
\hskip-3mm
\sum_{\lower5mm\hbox{$\scriptstyle
n^{(1)}\:{\rm in}\:(c_1,n_{1,0})\atop\scriptstyle\wt(n^{(1)})\,=\,w_1$}}
\hskip-6mm
q^{M+\sum_{i=1}^{c_1} n_{1,i}(n_{1,i}-1)}
\prod_{i=1}^{c_1}
\left[{2\left(w_1-\sum_{j=1}^i n_{1,j}\right)\!+n_{1,i-1}+n_{1,i}\atop
  n_{1,i-1}-n_{1,i}}\!\right]\!.
\end{eqnarray*}
These correspond to the terms for which $x=1$ in the premise.

On the other hand, for those cases where $n_{1,c_1}>0$,
use Lemma \ref{ShiftLemma}(4) and then the induction hypothesis:
\begin{eqnarray*}
\hbox to 5mm{$\displaystyle
R(N,M,0,p,1,p^\prime-1)
$\hss}\\[0.5mm]
&=&
\hskip-3mm
\sum_{\lower5mm\hbox{$\scriptstyle
n^{(1)}\:{\rm in}\:(c_1,n_{1,0})\atop\scriptstyle\wt(n^{(1)})\,\le\,w_1$}}
\hskip-6mm
\left(
q^{\sum_{i=1}^{c_1} n_{1,i}^2}
\prod_{i=1}^{c_1}
\left[{2\left(w_1-\sum_{j=1}^i n_{1,j}\right)+n_{1,i-1}+n_{1,i}\atop
  n_{1,i-1}-n_{1,i}}\right]
\right.\\[-2mm]
&&\left.\hskip10mm\phantom{M\atop N}\times
q^{w_1+1-n_{1,c_1}-\wt(n^{(1)})}
\right.\\[0.5mm]
&&\left.\hskip10mm\phantom{M\atop N}\times
R(w_1+1-\wt(n^{(1)}),n_{1,c_1}-1,0,p^\prime-c_1p,1,p-1)
\right)\\[0.5mm]
&=&
\hskip-3mm
\sum_{\lower5mm\hbox{$\scriptstyle
n^{(1)}\:{\rm in}\:(c_1,n_{1,0})\atop\scriptstyle\wt(n^{(1)})\,\le\,w_1$}}
\hskip-6mm
\left(
q^{\sum_{i=1}^{c_1} n_{1,i}^2}
\prod_{i=1}^{c_1}
\left[{2\left(w_1-\sum_{j=1}^i n_{1,j}\right)+n_{1,i-1}+n_{1,i}\atop
  n_{1,i-1}-n_{1,i}}\right]
\right.\\[-2mm]
&&\left.\hskip10mm\phantom{M\atop N}\times
q^{w_1+1-n_{1,c_1}-\wt(n^{(1)})+n_{1,c_1}-1}
\right.\\[0.5mm]
&&\left.\hskip10mm\phantom{M\atop N}\times
\hskip-31mm\sum_{\vbox{\openup1\jot\halign{\hfil
$\scriptstyle #$&$\scriptstyle #$\hfil\cr
\noalign{\vskip1mm}
&\!\!2\,\le\,x\,\le\,t+1\cr
n^{(2)}\:{\rm in}\:(c_2,n_{2,0}),\,&\wt(n^{(2)})\,\le\,w_2,\,n_{2,f_2}>0\cr
n^{(3)}\:{\rm in}\:(c_3,n_{3,0}),\,&\wt(n^{(3)})\,\le\,w_3,\,n_{3,f_3}>0\cr
&\:\vdots\cr
n^{(x-1)}\:{\rm in}\:(c_{x-1},n_{x-1,0}),\,&\wt(n^{(x-1)})\,\le\,w_{x-1},\,
n_{x-1,f_{x-1}}>0\cr
n^{(x)}\:{\rm in}\:(c_x,n_{x,0}),\,&\wt(n^{(x)})\,=\,w_x,\,
n_{x,f_x}=0\:{\rm if}\:x<t+1\cr}}}
\hskip-32mm\right.
\left(
\:q^{\sum_{k=2}^{x}\sum_{i=1}^{c_k} n_{k,i}(n_{k,i}-1)}
\phantom{M\atop N}\right.\\[0.5mm]
&&\left.\hskip15mm\phantom{M\atop N}\times
\prod_{k=2}^{x}\prod_{i=1}^{c_k}
\left[{2\left(w_k-\sum_{j=1}^i n_{k,j}\right)+n_{k,i-1}+n_{k,i}
\atop n_{k,i-1}-n_{k,i}}\right]
\right)\\[2mm]
&=&
\hskip7mm
q^{M}
\hskip-30mm\sum_{\vbox{\openup1\jot\halign{\hfil
$\scriptstyle #$&$\scriptstyle #$\hfil\cr
\noalign{\vskip1mm}
&\!\!2\,\le\,x\,\le\,t+1\cr
n^{(1)}\:{\rm in}\:(c_1,n_{1,0}),\,&\wt(n^{(1)})\,\le\,w_1,\,n_{1,f_1}>0\cr
n^{(2)}\:{\rm in}\:(c_2,n_{2,0}),\,&\wt(n^{(2)})\,\le\,w_2,\,n_{2,f_2}>0\cr
&\:\vdots\cr
n^{(x-1)}\:{\rm in}\:(c_{x-1},n_{x-1,0}),\,&\wt(n^{(x-1)})\,\le\,w_{x-1},\,
n_{x-1,f_{x-1}}>0\cr
n^{(x)}\:{\rm in}\:(c_x,n_{x,0}),\,&\wt(n^{(x)})\,=\,w_x,\,
n_{x,f_x}=0\:{\rm if}\:x<t+1\cr}}}
\hskip-32mm
\left(
\:q^{\sum_{k=1}^{x}\sum_{i=1}^{c_k} n_{k,i}(n_{k,i}-1)}
\phantom{M\atop N}\right.\\[0.5mm]
&&\left.\hskip15mm\phantom{M\atop N}\times
\prod_{k=1}^{x}\prod_{i=1}^{c_k}
\left[{2\left(w_k-\sum_{j=1}^i n_{k,j}\right)+n_{k,i-1}+n_{k,i}
\atop n_{k,i-1}-n_{k,i}}\right]
\right),
\end{eqnarray*}
after using $M=w_1$.
Combining this with the previous expression for when $n_{1,c_1}=0$,
proves the theorem in the case $t+1$.
Hence, by induction, the theorem holds for all $t\ge1$.
\cqfd

\begin{example}\label{ZeroEx} {\rm
We use the above theorem to calculate $R(6,5,0,6,1,16)$.
Via Lemma \ref{CFLemma}, the continued fraction of
$\cf{17}{6}$ is found to be $(2,1,5)$.
Thus $c_1=2$, $c_2=1$ and $c_3=3$.
For convenience, we consider separately the sequences of parameters
which correspond to the cases $x=1$, $x=2$ and $x=3$.
\begin{itemize}
\item $x=1$: require a 2-part partition $n^{(1)}=(n_{1,1},n_{1,2})$
having $n_{1,1}\le 6$, weight precisely 5 and $n_{1,2}=0$.
The only possibility is $n^{(1)}=(5,0)$.
Since $w_1=5$ and $n_{1,0}=6$, the corresponding summand of
Theorem \ref{DetGenThrm} is
$$
q^{20} \left[{11\atop1}\right] \left[{5\atop5}\right].
$$
\item $x=2$: require a 2-part partition $n^{(1)}=(n_{1,1},n_{1,2})$
satisfying $n_{1,1}\le 6$, $\wt(n^{(1)})\le 5$ and $n_{1,2}>0$, and
a 1-part partition $n^{(2)}=(n_{2,1})$ satisfying
$n_{2,1}\le 6-\wt(n^{(1)})$, $\wt(n^{(2)})=n_{1,2}-1$ and $n_{2,1}=0$.
The last two conditions on $n^{(2)}$ imply that $n_{1,2}=1$ and hence
the only possible sequences here are
$(n_{1,2},n_{1,2};n_{2,1})=(4,1;0)$, $(3,1;0)$, $(2,1;0)$ or
$(1,1;0)$.
As above, $w_1=5$ and $n_{1,0}=6$, and now $w_2=0$ in each case but
$n_{2,0}=6-\wt(n^{(1)})$.
Thereupon, the corresponding summands of Theorem \ref{DetGenThrm} are:
$$
q^{12} \left[{12\atop2}\right] \left[{5\atop3}\right] \left[{1\atop1}\right]
+q^{6} \left[{13\atop3}\right] \left[{6\atop2}\right] \left[{2\atop2}\right]
+q^{2} \left[{14\atop4}\right] \left[{7\atop1}\right] \left[{3\atop3}\right]
+\left[{15\atop5}\right] \left[{8\atop0}\right] \left[{4\atop4}\right].
$$
\item $x=3$: require a 2-part partition $n^{(1)}=(n_{1,1},n_{1,2})$
satisfying $n_{1,1}\le 6$, $\wt(n^{(1)})\le 5$ and $n_{1,2}>0$,
a 1-part partition $n^{(2)}=(n_{2,1})$ satisfying
$n_{2,1}\le 6-\wt(n^{(1)})$, $\wt(n^{(2)})\le n_{1,2}-1$ and $n_{2,1}>0$,
and a 3-part partition $n^{(3)}=(n_{3,1},n_{3,2},n_{3,3})$
satisfying $n_{3,1}\le n_{1,2}-\wt(n^{(2)})$ and $\wt(n^{(3)})=n_{2,1}-1$.
The only possible sequences here are
$(n_{1,2},n_{1,2};n_{2,1};$
$n_{3,1},n_{3,2},n_{3,3})=$
$(3,2;1;0,0,0)$ or $(2,2;1;0,0,0)$.
In each case, we have $w_1=5$, $w_2=1$, $w_3=0$, $n_{1,0}=6$ and
$n_{3,0}=1$. In the first case $n_{2,0}=1$ and in the second case
$n_{2,0}=2$.
Thereupon, the corresponding summands of Theorem \ref{DetGenThrm} are:
$$
q^{8} \left[{13\atop3}\right] \left[{5\atop1}\right] \left[{2\atop0}\right]
\left[{1\atop1}\right] \left[{0\atop0}\right] \left[{0\atop0}\right]
+q^{4} \left[{14\atop4}\right] \left[{6\atop0}\right] \left[{3\atop1}\right]
\left[{1\atop1}\right] \left[{0\atop0}\right] \left[{0\atop0}\right].
$$
\end{itemize}
Evaluating the Gaussians in all three cases, adding the results,
and multiplying by $q^M=q^5$ results in:
$$
q^5+q^6+3q^7+5q^8+10q^9+\cdots+1062q^{29}+1064q^{30}
+1062q^{31}+\cdots+3q^{53}+q^{54}+q^{55}.
$$
As may be verified, this agrees with the value obtained using
Theorem \ref{GenFunThrm}.
\endex}
\end{example}


\begin{theorem}\label{GenThrm}
Let $p$ and $p^\prime$ be positive coprime integers and let
$\cf{p^\prime}p$ have continued fraction
$(c_1,c_2,\ldots,c_{t-1},c_t+2)$ with $t\ge1$.
For $1\le u\le t$ and $1\le d\le c_{u}$,
let $r$ and $s$ be such that $(c_1,c_2,\ldots,c_{u-1},d)$ is the
continued fraction for $\cf sr$.
Then $R(N,M,r,p-r,s,p^\prime-s)$ may be expressed as follows:
\begin{eqnarray*}
\hbox to 5mm{$\displaystyle
R(N,M,r,p-r,s,p^\prime-s)
$\hss}\\[0.5mm]
&=&\sum
\left(
q^{n_{u,d} + \sum_{k=1}^{u} \sum_{i=1}^{f_k} n_{k,i}^2
+\sum_{k=u+1}^{x} \sum_{i=1}^{f_k} n_{k,i}(n_{k,i}-1)}
\phantom{M\atop N}\right. \\[0.5mm]
&&\left.\quad\qquad\qquad\qquad\times
\prod_{k=1}^{x}\prod_{i=1}^{f_k}
\left[{2\left(w_k-\sum_{j=1}^i n_{k,j}\right)+n_{k,i-1}+n_{k,i}
\atop n_{k,i-1}-n_{k,i}}\right]
\right),
\end{eqnarray*}
where, on setting
$(f_1,f_2,\ldots,f_t,f_{t+1})=
(c_1,c_2,\ldots,c_{u-1},d,c_{u}-d,c_{u+1},\ldots,c_t)$,
the sum is over all sequences
$n^{(1)}$, $n^{(2)},\ldots,n^{(x)}$ of partitions with $u\le x\le t+1$
for which,
for $1\le k\le x$,
the partition $n^{(k)}$ is in $(f_k,n_{k,0})$, satisfies
$\wt(n^{(k)})\le w_k$ with $\wt(n^{(x)})=w_x$,
and additionally also
satisfies $\wt(n^{(u)})<w_u$ if $u<x$,
satisfies $n_{k,f_k}>0$ if $u<k<x$,
and satisfies $n_{x,f_x}=0$ whenever $x<t+1$.
Here, we define:
$$
\mathalign{
w_1\:&=M;&\cr
w_k\:&=n_{k-1,f_{k-1}}&2\le k\le u;\cr
w_{u+1}\:&=w_u-1-\wt(n^{(u)});&\cr
w_k\:&=n_{k-1,f_{k-1}}-1\qquad\qquad &u+2\le k\le x;\cr}
$$
and:
$$
\mathalign{
n_{1,0}\:&=N;&\cr
n_{k,0}\:&=w_{k-1}-\wt(n^{(k-1)})
&2\le k\le u;\cr
n_{u+1,0}\:&=n_{u,f_u}+1;&\cr
n_{k,0}\:&=w_{k-1}+1-\wt(n^{(k-1)})\qquad
&u+2\le k\le x.\cr}
$$
\end{theorem}

\Proof
Let $\hat p$ and $\hat p^\prime$ be such that
$\cf{\hat p^\prime}{\hat p}$ has continued fraction
$(c_{u},c_{u+1},\ldots,c_{t}+2)$.
Then, on using $v=u$ and $d_v=d$ in Theorem \ref{CFThrm},
we obtain:
\begin{eqnarray*}
\hbox to 5mm{$\displaystyle
R(N,M,r,p-r,s,p^\prime-s)
$\hss}\\[0.5mm]
&=&
\sum
\left(\:
q^{\sum_{k=1}^{u-1}\sum_{i=1}^{c_k} n_{k,i}^2}
\prod_{k=1}^{u-1}\prod_{i=1}^{c_k}
\left[{2\left(w_k-\sum_{j=1}^i n_{k,j}\right)+n_{k,i-1}+n_{k,i}
\atop n_{k,i-1}-n_{k,i}}\right]
\right.\\[0.5mm]
&&\left.\hskip30mm\times\phantom{M\atop N}
R(n_{u,0},w_u,1,\hat p-1,d,\hat p^\prime-d)
\right),
\end{eqnarray*}
where the sum is over all sequences
$n^{(1)}$, $n^{(2)},\ldots,n^{(u-1)}$ of partitions for which, for $1\le k<u$,
the partition $n^{(k)}$ is in $(c_k,n_{k,0})$ and satisfies
$\wt(n^{(k)})\le w_k$.

Applying first Lemma \ref{OnePassLem} with $c=d=f_u$
to $R(n_{u,0},w_u,1,\hat p-1,d,\hat p^\prime-d)$,
and then Corollary \ref{GenFunCor} to the result, yields:
\begin{eqnarray*}
\hbox to 5mm{$\displaystyle
R(N,M,r,p-r,s,p^\prime-s)
$\hss}\\[0.5mm]
&=&
\sum
\left(\:
q^{\sum_{k=1}^{u}\sum_{i=1}^{f_k} n_{k,i}^2}
\prod_{k=1}^{u}\prod_{i=1}^{f_k}
\left[{2\left(w_k-\sum_{j=1}^i n_{k,j}\right)+n_{k,i-1}+n_{k,i}
\atop n_{k,i-1}-n_{k,i}}\right]
\right.\\[0.5mm]
&&\left.\hskip30mm\times\phantom{M\atop N}
R(n_{u,d},w_u-\wt(n^{(u)}),1,\hat p-1,0,\hat p^\prime-d\hat p)
\right)\\[0.5mm]
&=&
\sum
\left(\:
q^{\sum_{k=1}^{u}\sum_{i=1}^{f_k} n_{k,i}^2}
\prod_{k=1}^{u}\prod_{i=1}^{f_k}
\left[{2\left(w_k-\sum_{j=1}^i n_{k,j}\right)+n_{k,i-1}+n_{k,i}
\atop n_{k,i-1}-n_{k,i}}\right]
\right.\\[0.5mm]
&&\left.\hskip30mm\times\phantom{M\atop N}
R(w_u-\wt(n^{(u)}),n_{u,d},0,\hat p^\prime-d\hat p,1,\hat p-1)
\right),
\end{eqnarray*}
where, now, the sum is over all sequences
$n^{(1)}$, $n^{(2)},\ldots,n^{(u)}$ of partitions for which, for $1\le k\le u$,
the partition $n^{(k)}$ is in $(f_k,n_{k,0})$ and satisfies
$\wt(n^{(k)})\le w_k$.

Now, if $\wt(n^{(u)})=w_u$ then
$R(w_u-\wt(n^{(u)}),n_{u,d},0,\hat p^\prime-d\hat p,1,\hat p-1)
=\delta_{n_{u,d},0}$
by Lemma \ref{RootExcpt}.
Otherwise, we may transform this term using Lemma \ref{ShiftLemma}(4).
Thereupon,
\begin{eqnarray*}
\hbox to 1mm{$\displaystyle
R(N,M,r,p-r,s,p^\prime-s)
$\hss}\\[0.5mm]
&=&
\hskip-4mm
\sum_{\lower5mm\hbox{$\scriptstyle
\wt(n^{(u)})\,=\,w_u\atop\scriptstyle n_{u,d}=0$}}
\hskip-6mm
\left(\:
q^{\sum_{k=1}^{u}\sum_{i=1}^{f_k} n_{k,i}^2}
\prod_{k=1}^{u}\prod_{i=1}^{f_k}
\left[{2\left(w_k-\sum_{j=1}^i n_{k,j}\right)+n_{k,i-1}+n_{k,i}
\atop n_{k,i-1}-n_{k,i}}\right]
\right)\\[0.5mm]
&&\!\!+
\hskip-4mm
\sum_{\lower3mm\hbox{$\scriptstyle \wt(n^{(u)})\,<\,w_u$}}
\hskip-6mm
\left(\:
q^{\sum_{k=1}^{u}\sum_{i=1}^{f_k} n_{k,i}^2}
\prod_{k=1}^{u}\prod_{i=1}^{f_k}
\left[{2\left(w_k-\sum_{j=1}^i n_{k,j}\right)+n_{k,i-1}+n_{k,i}
\atop n_{k,i-1}-n_{k,i}}\right]
\right.\\[0.5mm]
&&\left.\hskip20mm\times\phantom{M\atop N}
q^{n_{u,d}+1-w_u+\wt(n^{(u)})}
\right.\\[0.5mm]
&&\left.\hskip20mm\times\phantom{M\atop N}
R(n_{u,d}+1,\,w_u-1-\wt(n^{(u)}),
0,\hat p,1,\hat p^\prime-d\hat p-1)\right),
\end{eqnarray*}
with the sums over all sequences as above.

Now, with $n_{u+1,0}=n_{u,d}+1$ and $w_{u+1}=w_u-1-\wt(n^{(u)})$,
we use the expression for
$R(n_{u+1,0},w_{u+1},0,\hat p,1,\hat p^\prime-d\hat p-1)$
given by Theorem \ref{DetGenThrm}.
First note that since $\cf{\hat p^\prime}{\hat p}$ has continued fraction
$(f_{u+1}+d,f_{u+2},\ldots,f_{t+1}+2)$, it follows that
the continued fraction of $\cf{\hat p^\prime-d\hat p}{\hat p}$
is $(f_{u+1},f_{u+2},\ldots,f_{t+1}+2)$.
Then Theorem \ref{DetGenThrm} yields:
\begin{eqnarray*}
\hbox to 5mm{$\displaystyle
R(n_{u+1,0},w_{u+1},0,\hat p,1,\hat p^\prime-d\hat p-1)
$\hss}\\[0.5mm]
&=&
\hskip7mm
q^{w_{u+1}}
\hskip-30mm\sum_{\vbox{\openup1\jot\halign{\hfil
$\scriptstyle #$&$\scriptstyle #$\hfil\cr
\noalign{\vskip1mm}
&\!\!u+1\,\le\,x\,\le\,t+1\cr
n^{(u+1)}\:{\rm in}\:(f_{u+1},n_{u+1,0}),\,&\wt(n^{(u+1)})\,\le\,w_{u+1},\,
n_{u+1,f_{u+1}}>0\cr
n^{(u+2)}\:{\rm in}\:(f_{u+2},n_{u+2,0}),\,&\wt(n^{(u+2)})\,\le\,w_{u+2},\,
n_{u+2,f_{u+2}}>0\cr
&\:\vdots\cr
n^{(x-1)}\:{\rm in}\:(f_{x-1},n_{x-1,0}),\,&\wt(n^{(x-1)})\,\le\,w_{x-1},\,
n_{x-1,f_{x-1}}>0\cr
n^{(x)}\:{\rm in}\:(f_x,n_{x,0}),\,&\wt(n^{(x)})\,=\,w_x,\,
n_{x,f_x}=0\:{\rm if}\:x<t+1\cr}}}
\hskip-32mm
\left(
\:q^{\sum_{k=u+1}^{x}\sum_{i=1}^{f_k} n_{k,i}(n_{k,i}-1)}
\phantom{M\atop N}\right.\\[0.5mm]
&&\left.\hskip15mm\phantom{M\atop N}\times
\prod_{k=u+1}^{x}\prod_{i=1}^{f_k}
\left[{2\left(w_k-\sum_{j=1}^i n_{k,j}\right)+n_{k,i-1}+n_{k,i}
\atop n_{k,i-1}-n_{k,i}}\right]
\right),
\end{eqnarray*}
with $n_{k,0}$ and $w_k$ for $u<k\le t+1$ as in the premise.
Combining this with the previous expression, proves the theorem.
\cqfd

\begin{example}\label{GenEx} {\rm
We use the above theorem to calculate $R(3,4,2,14,5,32)$.
Via Lemma \ref{CFLemma}, the continued fraction of
$\cf{37}{16}$ is found to be $(2,3,5)$, and the continued fraction
of $\cf52$ is found to be $(2,2)$. Thus here $t=3$, $u=2$, $d=2$,
$f_1=2$, $f_2=2$, $f_3=1$ and $f_4=3$.
Here, $x\ge2$ and so each partition sequence contains at least the two
partitions $n^{(1)}$ and $n^{(2)}$. These comprise $f_1=2$ and $f_2=2$
parts respectively, for which $n_{1,1}\le 3$ and $\wt(n^{(1)})\le4$,
$n_{2,1}\le 4-\wt(n^{(1)})$ and $\wt(n^{(2)})\le n_{1,2}$.
These constraints give rise to ten possible
$(n_{1,1},n_{1,2};n_{2,1},n_{2,2})$.
They are:
$$
\displaylines{
(3,0;0,0),\:
(2,1;1,0),\:
(2,0;0,0),\:
(1,1;1,0),\:
(1,0;0,0),\:
(0,0;0,0),
\cr
(3,1;0,0),\:
(2,1;0,0),\:
(1,1;0,0),\:
(2,2;0,0).
\cr}
$$
Since $w_u=n_{1,2}$, the first six cases have $\wt(n^{(u)})=w_u$
and the final four cases have $\wt(n^{(u)})>w_u$.
Thus the first six cases are precisely the partition sequences for $x=2$.

The remaining four cases give rise to those sequences for $x>2$.
First we require $n^{(3)}$ with one part $n_{3,1}$ for which
$n_{3,1}\le n_{3,0}$ where $n_{3,0}=n_{2,2}+1$ and
$\wt(n^{(3)})\le w_3$ where $w_3=n_{1,2}-1-\wt(n^{(2)})$.
In the first three cases, this gives $w_3=0$ and so $n_{3,1}=0$
in these cases. These are thus the cases for $x=3$.

For the final case, $n_{1,0}=1$ and $w_3=1$.
This cannot give rise to an $x=3$ term since the requirement that
$\wt(n^{(x)})=w_3$ forces $n_{3,1}=1$. For an $x=4$ term,
we also require $n_{3,1}=1$ here, whereupon $n_{4,0}=1$ and $w_4=0$,
so that only $n^{(4)}=(0,0,0)$ is possible.

The following table summarises these details.
\medskip
$$
\vbox{\offinterlineskip\halign{$\strut
\enspace #\hfil\enspace$&\vrule #&&$\hfil\enspace #\hfil$\cr
(n^{(1)};n^{(2)};\cdots;n^{(x)})
&&x&n_{1,0}&w_1&n_{2,0}&w_2&n_{3,0}&w_3&n_{4,0}&w_4\cr
\omit&height 1mm\cr
\multispan{11}\hrulefill\cr
\omit&height 1mm\cr
(3,0;0,0)&&2&3&4&1&0&&&&&\cr
(2,1;1,0)&&2&3&4&1&1\cr
(2,0;0,0)&&2&3&4&2&0\cr
(1,1;1,0)&&2&3&4&2&1\cr
(1,0;0,0)&&2&3&4&3&0\cr
(0,0;0,0)&&2&3&4&4&0\cr
(3,1;0,0;0)&&3&3&4&0&1&1&0\cr
(2,1;0,0;0)&&3&3&4&1&1&1&0\cr
(1,1;0,0;0)&&3&3&4&2&1&1&0\cr
(2,2;0,0;1;0,0,0)&&4&3&4&0&2&1&1&1&0&
\cr
}}
$$
\medskip
\noindent
Using this information to evaluate the expression in Theorem
\ref{GenThrm}, leads to:
$$
\displaylines{
q^9 \gauss80\gauss53\gauss11\gauss00
+q^6 \gauss91\gauss51\gauss20\gauss11
+q^4 \gauss91\gauss62\gauss22\gauss00
+q^3 \gauss{10}2\gauss60\gauss31\gauss11
\cr
+q \gauss{10}2\gauss71\gauss33\gauss00
+ \gauss{11}3\gauss80\gauss44\gauss00
+q^{10} \gauss80\gauss42\gauss20\gauss20\gauss11
+q^{5} \gauss91\gauss51\gauss31\gauss20\gauss11
\cr
+q^2 \gauss{10}2\gauss60\gauss42\gauss20\gauss11
+q^8 \gauss91\gauss40\gauss40\gauss40\gauss20\gauss11\gauss00\gauss00.
\cr}
$$
Evaluating all the Gaussians and summing, results in:
$$
1+2q+5q^2+10q^3+18q^4+28q^5+\cdots+
114q^{11}+119q^{12}+114q^{13}+\cdots+
2q^{23}+q^{24}.
$$
As may be verified, this agrees with the value obtained using
Theorem \ref{GenFunThrm}.
(In fact, from the viewpoint of partition pairs, this result is
fairly trivial --- being the generating function for {\em all}
pairs of partitions in $(3,4)\times(3,4)$.)
\endex}
\end{example}


\medskip\noindent
Note that, for convenience, the above theorem does not deal with the cases
of $(s:r)$ having continued fraction $()$ or
$(c_1,c_2,\ldots,c_{t-1},c_t+1)$. The former of these, when
$r=0$ and $s=1$, is already dealt with in Theorem \ref{DetGenThrm},
and as explained in Note \ref{Anote}, the latter is
equal to that resulting from the case $(c_1,c_2,\ldots,c_{t-1})$.

\section{Q-polynomial identities with a conventional bosonic side}

In the previous section, we obtained $q$-polynomial identities
that depend on two finite size parameters. We can take one of 
these parameters to infinity, and obtain identities that depend
on a single parameter. Depending on which parameter is taken to
infinity, we obtain polynomial identities that have a {\it conventional}
dependence on the rest of the parameters, namely $\{p, p', r, s\}$, 
or a {\it non-conventional} dependence. In this section, we present 
the former. In the next section, we present the latter.

To be more specific, we use the words {\it 'conventional'} and
{\it 'non-conventional'} in the following sense: We compare the 
dependence of the bosonic side of the $q$-polynomial identities 
that we obtain with that of the $q$-polynomial identities that
appear in {\it e.g.} \cite{bms}. 


\begin{theorem}\label{NinfThrm}
Let $p$ and $p^\prime$ be positive coprime integers and let
$\cf{p^\prime}p$ have continued fraction
$(c_1,c_2,\ldots,c_{t-1},c_t+2)$ with $t\ge1$.
For $1\le u\le t$ and $1\le d\le c_{u}$,
let $r$ and $s$ be such that $(c_1,c_2,\ldots,c_{u-1},d)$ is the
continued fraction for $\cf sr$.
Then
\begin{eqnarray*}
\hbox to 5mm{$\displaystyle
\lim_{N\to\infty}R(N,M,r,p-r,s,p^\prime-s)
$\hss}\\[0.5mm]
&=&{1\over (q)_{2M}}\sum
\left(
q^{n_{u,d} + \sum_{k=1}^{u} \sum_{i=1}^{f_k} n_{k,i}^2
+\sum_{k=u+1}^{x} \sum_{i=1}^{f_k} n_{k,i}(n_{k,i}-1)}
\phantom{M\atop N}\right. \\[0.5mm]
&&\left.\qquad\qquad\times
\prod_{k=1}^{x}\prod_{i=1+\delta_{1,k}}^{f_k}
\left[{2\left(w_k-\sum_{j=1}^i n_{k,j}\right)+n_{k,i-1}+n_{k,i}
\atop n_{k,i-1}-n_{k,i}}\right]
\right),
\end{eqnarray*}
where, on setting
$(f_1,f_2,\ldots,f_t,f_{t+1})=
(c_1,c_2,\ldots,c_{u-1},d,c_{u}-d,c_{u+1},\ldots,c_t)$,
the sum is over all sequences
$n^{(1)}$, $n^{(2)},\ldots,n^{(x)}$ of partitions with $u\le x\le t+1$
for which,
for $1\le k\le x$,
the partition $n^{(k)}$ is in $(f_k,n_{k,0})$, satisfies
$\wt(n^{(k)})\le w_k$ with $\wt(n^{(x)})=w_x$,
and additionally also
satisfies $\wt(n^{(u)})<w_u$ if $u<x$,
satisfies $n_{k,f_k}>0$ if $u<k<x$,
and satisfies $n_{x,f_x}=0$ whenever $x<t+1$.
Here, we define:
$$
\mathalign{
w_1\:&=M;&\cr
w_k\:&=n_{k-1,f_{k-1}}&2\le k\le u;\cr
w_{u+1}\:&=w_u-1-\wt(n^{(u)});&\cr
w_k\:&=n_{k-1,f_{k-1}}-1\qquad\qquad &u+2\le k\le x;\cr}
$$
and:
$$
\mathalign{
n_{1,0}\:&=\infty;&\cr
n_{k,0}\:&=w_{k-1}-\wt(n^{(k-1)})
&2\le k\le u;\cr
n_{u+1,0}\:&=n_{u,f_u}+1;&\cr
n_{k,0}\:&=w_{k-1}+1-\wt(n^{(k-1)})\qquad
&u+2\le k\le x.\cr}
$$
\end{theorem}

\Proof Only the $i=k=1$ Gaussian in Theorem \ref{GenThrm} involves
$N$ ($N=n_{1,0}$).
On using $M=n_{0,f_0}$,
this term is
$$
\left[{2M + N - n_{1,1}\atop N-n_{1,1}}\right]
={(q)_{2M+N-n_{1,1}}\over (q)_{2M} (q)_{N-n_{1,1}}}.
$$
In the limit as $N\to\infty$, this tends termwise to $1/(q)_{2M}$.
The theorem then follows.
\cqfd

\begin{example}\label{NinfEx} {\rm
To illustrate the above theorem, we extend Example \ref{GenEx}
to calculate $\lim_{N\to\infty}R(N,4,2,14,5,32)$.
As in Example \ref{GenEx}, the continued fraction of
$\cf{37}{16}$ is $(2,3,5)$, the continued fraction
of $\cf52$ is $(2,2)$, and
$t=3$, $u=2$, $d=2$, $f_1=2$, $f_2=2$, $f_3=1$ and $f_4=3$.
Now, although the largest part of $n^{(1)}$ is unbounded,
the weight of $n^{(1)}$ is bounded by $w_1=4$ as before.
Thus, compared with Example \ref{GenEx}, only one more
partition sequence is permitted: it is $(4,0;0,0)$.
The full set of valid partition sequences are then summarised
in the following table.
\medskip
$$
\vbox{\offinterlineskip\halign{$\strut
\enspace #\hfil\enspace$&\vrule #&&$\hfil\enspace #\hfil$\cr
(n^{(1)};n^{(2)};\cdots;n^{(x)})
&&x&n_{1,0}&w_1&n_{2,0}&w_2&n_{3,0}&w_3&n_{4,0}&w_4\cr
\omit&height 1mm\cr
\multispan{11}\hrulefill\cr
\omit&height 1mm\cr
(4,0;0,0)&&2&\infty&4&0&0\cr
(3,0;0,0)&&2&\infty&4&1&0\cr
(2,1;1,0)&&2&\infty&4&1&1\cr
(2,0;0,0)&&2&\infty&4&2&0\cr
(1,1;1,0)&&2&\infty&4&2&1\cr
(1,0;0,0)&&2&\infty&4&3&0\cr
(0,0;0,0)&&2&\infty&4&4&0\cr
(3,1;0,0;0)&&3&\infty&4&0&1&1&0\cr
(2,1;0,0;0)&&3&\infty&4&1&1&1&0\cr
(1,1;0,0;0)&&3&\infty&4&2&1&1&0\cr
(2,2;0,0;1;0,0,0)&&4&\infty&4&0&2&1&1&1&0\cr
}}
$$
\medskip
\noindent
The calculation of $\,\lim_{N\to\infty}R(N,4,2,14,5,32)$
using Theorem \ref{NinfThrm} is then similar in detail to that
of Example \ref{GenEx}. We obtain:
\begin{eqnarray*}
\hbox to 0mm{$\displaystyle
\lim_{N\to\infty}R(N,4,2,14,5,32)
$\hss}\\[0.5mm]
&=&
{1\over(q)_8}\left(
q^{16} \gauss44\gauss00\gauss00
+q^9 \gauss53\gauss11\gauss00
+q^6 \gauss51\gauss20\gauss11
+q^4 \gauss62\gauss22\gauss00
\right.
\\[0.5mm]
&&
\:{} +q^3 \gauss60\gauss31\gauss11
+q \gauss71\gauss33\gauss00
+ \gauss80\gauss44\gauss00
+q^{10} \gauss42\gauss20\gauss20\gauss11
\\[0.5mm]
&&
\left.
\:{} +q^{5} \gauss51\gauss31\gauss20\gauss11
+q^2 \gauss60\gauss42\gauss20\gauss11
+q^8 \gauss40\gauss40\gauss40\gauss20\gauss11\gauss00\gauss00
\right).
\\[1mm]
&=&
(1+q+2q^2+3q^3+5q^4+5q^5+7q^6+7q^7+8q^8+\cdots+q^{15}+q^{16})
/(q)_8
\\[0.5mm]
&=&
1+2q+5q^2+10q^3+20q^4+34q^5+59q^6+94q^7+149q^8+224q^9+\cdots
\end{eqnarray*}
\endex}
\end{example}

\noindent From the above Theorem \ref{NinfThrm} together with
Theorem \ref{GenFunThrm}, one obtains the following
conventional $q$-polynomial identity 
for a finite version of the character $\chi^{p, p'}_{r, s}$, 
where the finite size parameter is $N$.

\begin{corollary}\label{FirstCor}
Let $\{p, p', r, s\}$ be as in Theorem \ref{NinfThrm}. Then
\begin{eqnarray*}
\hbox to 5mm{$\displaystyle
\sum_{k=-\infty}^\infty
q^{k^2pp'+k(p'r-ps)}
\left[{2 M \atop M + p'k}\right]
-
q^{(kp+r)(kp'+s)}
\left[{2 M \atop M + p' k + s }\right]
$\hss}
\\[1mm]
&&=\quad
\sum
\left(
q^{n_{u,d} + \sum_{k=1}^{u} \sum_{i=1}^{f_k} n_{k,i}^2
+\sum_{k=u+1}^{x} \sum_{i=1}^{f_k} n_{k,i}(n_{k,i}-1)}
\phantom{M\atop N}\right. \\[-0.5mm]
&&\left.\qquad\qquad\times\:
\prod_{k=1}^{x}\prod_{i=1+\delta_{1,k}}^{f_k}
\left[{2\left(w_k-\sum_{j=1}^i n_{k,j}\right)+n_{k,i-1}+n_{k,i}
\atop n_{k,i-1}-n_{k,i}}\right]
\right),
\end{eqnarray*}
where the sum is over the same parameters as that in Theorem \ref{NinfThrm}. 
\end{corollary}

\begin{note}{\rm Notice the dependence of the Gaussian polynomials, on 
the bosonic side of the above identities, on the parameters 
$\{p, p', r, s\}$: only $\{p', s\}$ appear, but not $\{p, r\}$. This 
is what we refer to as a {\it conventional} dependence on the parameters. 
Also notice that the form of the fermionic side is the same as appears 
in {\it e.g.} \cite{bms}.}
\end{note}

\section{Q-polynomial identities with a non-conventional bosonic side}

In this section, we take the finite size parameter, that was left
fixed in the previous section, to infinity, and retain the other 
fixed. We obtain another set of $q$-polynomial identities, that 
depend on a single finite size parameter, and that can be characterised 
by the fact that their bosonic side has a non-conventional dependence 
on the parameters $\{p, p', r, s\}$ that characterise the model. 

\begin{theorem}\label{MinfThrm}
Let $p$ and $p^\prime$ be positive coprime integers and let
$\cf{p^\prime}p$ have continued fraction
$(c_1,c_2,\ldots,c_{t-1},c_t+2)$ with $t\ge1$.
For $1\le u\le t$ and $1\le d\le c_{u}$,
let $r$ and $s$ be such that $(c_1,c_2,\ldots,c_{u-1},d)$ is the
continued fraction for $\cf sr$.
Then if $u>1$,
\begin{eqnarray*}
\hbox to 5mm{$\displaystyle
\lim_{M\to\infty}R(N,M,r,p-r,s,p^\prime-s)
$\hss}\\[0.5mm]
&=&\sum
\left(
q^{n_{u,d} + \sum_{k=1}^{u} \sum_{i=1}^{f_k} n_{k,i}^2
+\sum_{k=u+1}^{x} \sum_{i=1}^{f_k} n_{k,i}(n_{k,i}-1)}
\phantom{M\atop N}\right. \\[0.5mm]
&&\left.\qquad\times\phantom{M\atop N}
{1\over (q)_{2n_{1,f_1}}}
\prod_{i=1}^{f_1} {1\over (q)_{n_{1,i-1}-n_{1,i}}}
\right. \\[0.5mm]
&&\left.\qquad\times\phantom{M\atop N}
\prod_{k=2}^{x}\prod_{i=1+\delta_{2,k}}^{f_k}
\left[{2\left(w_k-\sum_{j=1}^i n_{k,j}\right)+n_{k,i-1}+n_{k,i}
\atop n_{k,i-1}-n_{k,i}}\right]
\right),
\end{eqnarray*}
otherwise if $u=1$,
\begin{eqnarray*}
\hbox to 5mm{$\displaystyle
\lim_{M\to\infty}R(N,M,r,p-r,s,p^\prime-s)
$\hss}\\[0.5mm]
&=&\sum
\left(
q^{n_{1,d} + \sum_{i=1}^{f_1} n_{1,i}^2
+\sum_{k=2}^{x} \sum_{i=1}^{f_k} n_{k,i}(n_{k,i}-1)}
\phantom{M\atop N}\right. \\[0.5mm]
&&\left.\qquad\times\phantom{M\atop N}
{1\over (q)_{2n_{2,f_2}}}
\prod_{i=1}^{f_1} {1\over (q)_{n_{1,i-1}-n_{1,i}}}
\prod_{i=1}^{f_2} {1\over (q)_{n_{2,i-1}-n_{2,i}}}
\right. \\[0.5mm]
&&\left.\qquad\times\phantom{M\atop N}
\prod_{k=3}^{x}\prod_{i=1+\delta_{3,k}}^{f_k}
\left[{2\left(w_k-\sum_{j=1}^i n_{k,j}\right)+n_{k,i-1}+n_{k,i}
\atop n_{k,i-1}-n_{k,i}}\right]
\right),
\end{eqnarray*}
where, on setting
$(f_1,f_2,\ldots,f_t,f_{t+1})=
(c_1,c_2,\ldots,c_{u-1},d,c_{u}-d,c_{u+1},\ldots,c_t)$,
the sum is over all sequences
$n^{(1)}$, $n^{(2)},\ldots,n^{(x)}$ of partitions with $u\le x\le t+1$
for which,
for $1\le k\le x$,
the partition $n^{(k)}$ is in $(f_k,n_{k,0})$, satisfies
$\wt(n^{(k)})\le w_k$ with $\wt(n^{(x)})=w_x$,
and additionally also
satisfies $\wt(n^{(u)})<w_u$ if $u<x$,
satisfies $n_{k,f_k}>0$ if $u<k<x$,
and satisfies $n_{x,f_x}=0$ whenever $x<t+1$.
Here, we define:
$$
\mathalign{
w_1\:&=\infty;&\cr
w_k\:&=n_{k-1,f_{k-1}}&2\le k\le u;\cr
w_{u+1}\:&=w_u-1-\wt(n^{(u)});&\cr
w_k\:&=n_{k-1,f_{k-1}}-1\qquad\qquad &u+2\le k\le x;\cr}
$$
and:
$$
\mathalign{
n_{1,0}\:&=N;&\cr
n_{k,0}\:&=w_{k-1}-\wt(n^{(k-1)})
&2\le k\le u;\cr
n_{u+1,0}\:&=n_{u,f_u}+1;&\cr
n_{k,0}\:&=w_{k-1}+1-\wt(n^{(k-1)})\qquad
&u+2\le k\le x.\cr}
$$
\end{theorem}

\Proof Consider the expression given by Theorem \ref{GenThrm}.
In the case $u>1$, the only Gaussians involving $M$ are those
factors indexed by $k=1$ and $i$ for $1\le i\le f_1$
(since $w_1=M$), and that
indexed by $k=2$ and $i=1$ (since $n_{2,0}=M+\sum_{i=1}^{f_1}n_{1,i}$).
All but the last of these terms is:
$$
\left[{2M + n_{1,i-1} + n_{1,i} - 2\sum_{j=1}^{i} n_{1,i}
\atop n_{1,i-1}-n_{1,i}}\right]
={(q)_{2M + n_{1,i-1} + n_{1,i} - 2\sum_{j=1}^{i} n_{1,i}}
\over (q)_{n_{1,i-1}-n_{1,i}} (q)_{2M-2\sum_{j=1}^{i-1} n_{1,i}}}.
$$
In the limit as $M\to\infty$, these tend termwise to
$1/(q)_{n_{1,i-1}-n_{1,i}}$.
The $k=2$, $i=1$ term is:
$$
\left[{2w_2+n_{2,0}-n_{2,1} \atop n_{2,0}-n_{2,1}}\right]
={(q)_{2w_2+n_{2,0}-n_{2,1}}
\over (q)_{n_{2,0}-n_{2,1}} (q)_{2w_2}}.
$$
As $M\to\infty$, then $n_{2,0}\to\infty$ whereupon this term
tends termwise to
$1/(q)_{2w_2}=1/(q)_{2n_{1,f_1}}$.

In the case where $u=1$, we have $w_1=M$, $w_2=M-\sum_{i=1}^{f_1} n_{1,i}$,
$n_{2,0}=n_{1,f_1}$ and
$n_{3,0}=M-1-\sum_{i=1}^{f_1} n_{1,i}-\sum_{i=1}^{f_2} n_{2,i}$.
Thus as $M\to\infty$, all the $k=1$ and $k=2$ terms behave similar
to the $k=1$ terms above,
and the $k=3$, $i=1$ term behaves as did the $k=2$, $i=1$ term above.
The theorem then follows.
\cqfd

\begin{example}\label{MinfEx} {\rm
To illustrate the above theorem, we extend Example \ref{GenEx}
to calculate $\lim_{M\to\infty}R(3,M,2,14,5,32)$.
As in Example \ref{GenEx}, the continued fraction of
$\cf{37}{16}$ is $(2,3,5)$, the continued fraction
of $\cf52$ is $(2,2)$, and
$t=3$, $u=2$, $d=2$, $f_1=2$, $f_2=2$, $f_3=1$ and $f_4=3$.
Now the weight of $n^{(1)}$ is unbounded. However,
the largest part of $n^{(1)}$ is bounded by $n_{1,0}=3$,
thus restricting the possible $n^{(1)}$ to a finite number.
In addition, $n_{2,0}=\infty$ and $w_2$ is finite.
Thus the weight of $n^{(2)}$ is bounded, thereby restricting the
possible $n^{(2)}$ to a finite number.
Thus overall, the permitted number of partition sequences is finite.
In fact, as may be readily confirmed, there are 21 of them.
They are tabulated in Appendix B.
Thereupon, as shown in Appendix B, Theorem \ref{MinfThrm} yields:
\begin{eqnarray*}
\hbox to 0mm{$\displaystyle
\lim_{M\to\infty}R(3,M,2,14,5,32)
$\hss}\\[0.5mm]
&=&
1+2q+5q^2+10q^3+18q^4+30q^5+49q^6+74q^7+110q^8+158q^9+\cdots
\end{eqnarray*}
\endex}
\end{example}

\noindent From the above Theorem \ref{MinfThrm} together with
Theorem \ref{GenFunThrm}, one obtains the following
$q$-polynomial identity 
for a finite version of the character $\chi^{p, p'}_{r, s}$, 
where the finite size parameter is $M$.
In this case, the expression so obtained has both the bosonic
and fermionic sides in a non-conventional form.

\begin{corollary}\label{SecondCor}
Let $\{p, p', r, s\}$ be as in Theorem \ref{MinfThrm}. 
If $u>1$ (as defined in Theorem \ref{MinfThrm}), then
\begin{eqnarray*}
\hbox to 5mm{$\displaystyle
\sum_{k=-\infty}^\infty
q^{k^2pp'+k(p'r-ps)}
\left[{2 N \atop N + p k}\right]
-
q^{(kp+r)(kp'+s)}
\left[{ 2 N \atop N + p k + r}\right]
$\hss}\\[1mm]
&&=
\quad
(q)_{2N}
\sum\:
\left(
q^{n_{u,d} + \sum_{k=1}^{u} \sum_{i=1}^{f_k} n_{k,i}^2
+\sum_{k=u+1}^{x} \sum_{i=1}^{f_k} n_{k,i}(n_{k,i}-1)}
\phantom{M\atop N}\right. \\[0.5mm]
&&\left.\qquad\qquad\times\phantom{M\atop N}
{1\over (q)_{2n_{1,f_1}}}
\prod_{i=1}^{f_1} {1\over (q)_{n_{1,i-1}-n_{1,i}}}
\right. \\[0.5mm]
&&\left.\qquad\qquad\times\phantom{M\atop N}
\prod_{k=2}^{x}\prod_{i=1+\delta_{2,k}}^{f_k}
\left[{2\left(w_k-\sum_{j=1}^i n_{k,j}\right)+n_{k,i-1}+n_{k,i}
\atop n_{k,i-1}-n_{k,i}}\right]
\right),
\end{eqnarray*}
otherwise if $u=1$,
\begin{eqnarray*}
\hbox to 5mm{$\displaystyle
\sum_{k=-\infty}^\infty
q^{k^2pp'+k(p'r-ps)}
\left[{2 N \atop N + p k}\right]
-
q^{(kp+r)(kp'+s)}
\left[{ 2 N \atop N + p k + r}\right]
$\hss}\\[1mm]
&&=
\quad
(q)_{2N}
\sum\:
\left(
q^{n_{1,d} + \sum_{i=1}^{f_1} n_{1,i}^2
+\sum_{k=2}^{x} \sum_{i=1}^{f_k} n_{k,i}(n_{k,i}-1)}
\phantom{M\atop N}\right. \\[0.5mm]
&&\left.\qquad\qquad\times\phantom{M\atop N}
{1\over (q)_{2n_{2,f_2}}}
\prod_{i=1}^{f_1} {1\over (q)_{n_{1,i-1}-n_{1,i}}}
\prod_{i=1}^{f_2} {1\over (q)_{n_{2,i-1}-n_{2,i}}}
\right. \\[0.5mm]
&&\left.\qquad\qquad\times\phantom{M\atop N}
\prod_{k=3}^{x}\prod_{i=1+\delta_{3,k}}^{f_k}
\left[{2\left(w_k-\sum_{j=1}^i n_{k,j}\right)+n_{k,i-1}+n_{k,i}
\atop n_{k,i-1}-n_{k,i}}\right]
\right),
\end{eqnarray*}
where the sums are taken over the same parameters as in 
Theorem \ref{MinfThrm}.
\end{corollary}

\begin{note}{\rm Notice the dependence of the Gaussian 
polynomials, on the bosonic side of the above identities, on the 
parameters $\{p, p', r, s\}$: only $\{p, r\}$ appear, but not 
$\{p', s\}$. This is what we refer to as a {\it non-conventional} 
dependence on the parameters. Also notice the form of the fermionic 
side: It is {\it not} the same form as appears in {\it e.g.} \cite{bms}.
In fact, it is reminiscent of what results from the application of 
the Bailey transform, see {\it e.g.} \cite{fq}.}
\end{note}

\section{Q-series identities}

In this section, we take both finite size parameters to infinity,
and obtain $q$-series identities. Notice that, we obtain the same
result, irrespectively of the order of removing the finite size
parameters.

\begin{theorem}\label{DinfThrm}
Let $p$ and $p^\prime$ be positive coprime integers and let
$\cf{p^\prime}p$ have continued fraction
$(c_1,c_2,\ldots,c_{t-1},c_t+2)$ with $t\ge1$.
For $1\le u\le t$ and $1\le d\le c_{u}$,
let $r$ and $s$ be such that $(c_1,c_2,\ldots,c_{u-1},d)$ is the
continued fraction for $\cf sr$.
Then if $u>1$,
\begin{eqnarray*}
\hbox to 5mm{$\displaystyle
\lim_{\scriptstyle M\to\infty\atop \scriptstyle N\to\infty}
R(N,M,r,p-r,s,p^\prime-s)
$\hss}\\[0.5mm]
&=&{1\over (q)_\infty}\sum
\left(
q^{n_{u,d} + \sum_{k=1}^{u} \sum_{i=1}^{f_k} n_{k,i}^2
+\sum_{k=u+1}^{x} \sum_{i=1}^{f_k} n_{k,i}(n_{k,i}-1)}
\phantom{M\atop N}\right. \\[0.5mm]
&&\left.\qquad\times\phantom{M\atop N}
{1\over (q)_{2n_{1,f_1}}}
\prod_{i=2}^{f_1} {1\over (q)_{n_{1,i-1}-n_{1,i}}}
\right. \\[0.5mm]
&&\left.\qquad\times\phantom{M\atop N}
\prod_{k=2}^{x}\prod_{i=1+\delta_{2,k}}^{f_k}
\left[{2\left(w_k-\sum_{j=1}^i n_{k,j}\right)+n_{k,i-1}+n_{k,i}
\atop n_{k,i-1}-n_{k,i}}\right]
\right),
\end{eqnarray*}
otherwise if $u=1$,
\begin{eqnarray*}
\hbox to 5mm{$\displaystyle
\lim_{\scriptstyle M\to\infty\atop \scriptstyle N\to\infty}
R(N,M,r,p-r,s,p^\prime-s)
$\hss}\\[0.5mm]
&=&{1\over (q)_\infty}\sum
\left(
q^{n_{1,d} + \sum_{i=1}^{f_1} n_{1,i}^2
+\sum_{k=2}^{x} \sum_{i=1}^{f_k} n_{k,i}(n_{k,i}-1)}
\phantom{M\atop N}\right. \\[0.5mm]
&&\left.\qquad\times\phantom{M\atop N}
{1\over (q)_{2n_{2,f_2}}}
\prod_{i=2}^{f_1} {1\over (q)_{n_{1,i-1}-n_{1,i}}}
\prod_{i=1}^{f_2} {1\over (q)_{n_{2,i-1}-n_{2,i}}}
\right. \\[0.5mm]
&&\left.\qquad\times\phantom{M\atop N}
\prod_{k=3}^{x}\prod_{i=1+\delta_{3,k}}^{f_k}
\left[{2\left(w_k-\sum_{j=1}^i n_{k,j}\right)+n_{k,i-1}+n_{k,i}
\atop n_{k,i-1}-n_{k,i}}\right]
\right),
\end{eqnarray*}
where, on setting
$(f_1,f_2,\ldots,f_t,f_{t+1})=
(c_1,c_2,\ldots,c_{u-1},d,c_{u}-d,c_{u+1},\ldots,c_t)$,
the sum is over all sequences
$n^{(1)}$, $n^{(2)},\ldots,n^{(x)}$ of partitions with $u\le x\le t+1$
for which,
for $1\le k\le x$,
the partition $n^{(k)}$ is in $(f_k,n_{k,0})$, satisfies
$\wt(n^{(k)})\le w_k$ with $\wt(n^{(x)})=w_x$,
and additionally also
satisfies $\wt(n^{(u)})<w_u$ if $u<x$,
satisfies $n_{k,f_k}>0$ if $u<k<x$,
and satisfies $n_{x,f_x}=0$ whenever $x<t+1$.
Here, we define:
$$
\mathalign{
w_1\:&=\infty;&\cr
w_k\:&=n_{k-1,f_{k-1}}&2\le k\le u;\cr
w_{u+1}\:&=w_u-1-\wt(n^{(u)});&\cr
w_k\:&=n_{k-1,f_{k-1}}-1\qquad\qquad &u+2\le k\le x;\cr}
$$
and:
$$
\mathalign{
n_{1,0}\:&=\infty;&\cr
n_{k,0}\:&=w_{k-1}-\wt(n^{(k-1)})
&2\le k\le u;\cr
n_{u+1,0}\:&=n_{u,f_u}+1;&\cr
n_{k,0}\:&=w_{k-1}+1-\wt(n^{(k-1)})\qquad
&u+2\le k\le x.\cr}
$$
\end{theorem}

\Proof In the expression given by Theorem \ref{MinfThrm}
(there is no difficulty in taking either limit first),
only the term $1/(q)_{n_{1,0}-n_{1,1}}=1/(q)_{N-n_{1,1}}$
is affected on taking $N\to\infty$. This gives rise to
$1/(q)_\infty$, with all other terms being unaffected.
\cqfd

\begin{example}\label{DinfEx} {\rm
To illustrate Theorem \ref{DinfThrm}, we consider, as in
Examples \ref{GenEx}, \ref{NinfEx} and \ref{MinfEx},
the case of $p=16$, $p^\prime=37$, $r=2$ and $s=5$.
As in those previous Examples, the continued fraction of
$\cf{37}{16}$ is $(2,3,5)$, the continued fraction
of $\cf52$ is $(2,2)$,
$t=3$, $u=2$, $d=2$, $f_1=2$, $f_2=2$, $f_3=1$ and $f_4=3$.
Unfortunately, unlike the previous examples, the number
of partition sequences over which the sum is taken is not finite,
since both $n_{1,0}=\infty$ and $w_1=\infty$.
Thus $n^{(1)}$ is any two-part partition $n^{(1)}=(n_{1,1},n_{1,2})$.
However, only a finite number of terms are required to guarantee
accuracy to any particular order. (In view of the exponent of
$q$ in the expression in Theorem \ref{DinfThrm}, if accuracy is
required up to the term of order $z$, then only those partitions
for which $\sum_{i=1}^{f_1} n_{1,i}^2 \le z$ are required.)
In the current case, accuracy up to the term of order 10 is guaranteed
by taking those $n_{1,1}$ and $n_{1,2}$ lying in the plane,
on or within the circle of radius $\sqrt{10}$, and in the sector for
which $n_{1,1}\ge n_{1,2}\ge0$.
There are eight such partitions.
They are:
$$
(3,1),\:
(3,0),\:
(2,2),\:
(2,1),\:
(2,0),\:
(1,1),\:
(1,0),\:
(0,0).
$$
Then, since $w_2=n_{1,2}$, the number of possible $n^{(2)}$ is
bounded in each case. Altogether, there are thirteen possible
partition sequences with $n^{(1)}$ as one of the above.
They are the sequences marked with an asterisk in the table
in Appendix B.
Using these in the expression in Theorem \ref{DinfThrm} yields:
\begin{eqnarray*}
\hbox to 0mm{$\displaystyle
\lim_{\scriptstyle M\to\infty\atop \scriptstyle N\to\infty}
R(N,M,2,14,5,32)
$\hss}\\[0.5mm]
&=&
\!{1\over (q)_\infty}
\left(
q^{11} \qeta2 \qeta2 \gauss11
+ q^9    \qeta0 \qeta3 \gauss00
+ q^{12} \qeta4 \qeta0 \gauss22
+ q^6    \qeta2 \qeta1 \gauss11
\right.
\\[0.5mm]
&&
\qquad\quad{}+ q^4    \qeta0 \qeta2 \gauss00
+ q^3    \qeta2 \qeta0 \gauss11
+ q      \qeta0 \qeta1 \gauss00
+        \qeta0 \qeta0 \gauss00
\\[1mm]
&&
\qquad\quad{}+ q^{10}    \qeta2 \qeta2 \gauss20 \gauss11
+ q^9    \qeta4 \qeta0 \gauss31 \gauss11
+ q^5    \qeta2 \qeta1 \gauss20 \gauss11
\\[0.5mm]
&&
\left.
\qquad\quad{}+ q^2    \qeta2 \qeta0 \gauss20 \gauss11
+ q^8    \qeta4 \qeta0 \gauss40 \gauss20 \gauss11 \gauss00 \gauss00
\right)
\\[0.5mm]
&=&
\!{1\over (q)_\infty}
\left(1\!+\!q\!+\!2q^2\!+\!3q^3\!+\!5q^4\!+\!7q^5\!+\!11q^6\!+\!15q^7
\!+\!22q^8\!+\!30q^9\!+\!41q^{10}\!+\cdots\:\right)
\\[0.5mm]
&=&
\!1\!+\!2q\!+\!5q^2\!+\!10q^3\!+\!20q^4\!+\!36q^5\!+\!65q^6\!+\!110q^7
\!+\!185q^8\!+\!300q^9\!+\!480q^{10}\!+\cdots
\end{eqnarray*}
Note that this generating function differs from that of the
unrestricted partition pairs, first at the $q^{10}$ term.
This reflects the fact that the
smallest partition pair that does not satisfy the constraints
of this example is $((5,5),(0,0))$.
\endex}
\end{example}

\noindent The above Theorem \ref{DinfThrm} together with
(\ref{RochaEq}) and (\ref{BRoachEq})
results in the following $q$-series identities
for the Virasoro character $\chi^{p, p'}_{r, s}$.

\begin{corollary}\label{ThirdCor}
Let $\{p, p', r, s\}$ be as in Theorem \ref{DinfThrm}. If $u>1$ 
(as defined in Theorem \ref{DinfThrm}), then
\begin{eqnarray*}
{1\over (q)_\infty}
&\phantom{-}&\hskip-9mm
\left(
\sum_{k=-\infty}^\infty
q^{k^2pp'+k(p'r-ps)}
-
q^{(kp+r)(kp'+s)}\right)\\
&=&
\sum
\left(
q^{n_{u,d} + \sum_{k=1}^{u} \sum_{i=1}^{f_k} n_{k,i}^2
+\sum_{k=u+1}^{x} \sum_{i=1}^{f_k} n_{k,i}(n_{k,i}-1)}
\phantom{M\atop N}\right. \\[0.5mm]
&&\left.\qquad\times\phantom{M\atop N}
{1\over (q)_{2n_{1,f_1}}}
\prod_{i=2}^{f_1} {1\over (q)_{n_{1,i-1}-n_{1,i}}}
\right. \\[0.5mm]
&&\left.\qquad\times\phantom{M\atop N}
\prod_{k=2}^{x}\prod_{i=1+\delta_{2,k}}^{f_k}
\left[{2\left(w_k-\sum_{j=1}^i n_{k,j}\right)+n_{k,i-1}+n_{k,i}
\atop n_{k,i-1}-n_{k,i}}\right]
\right),
\end{eqnarray*}
otherwise if $u=1$,
\begin{eqnarray*}
{1\over (q)_\infty}
&\phantom{-}&\hskip-9mm
\left(
\sum_{k=-\infty}^\infty
q^{k^2pp'+k(p'r-ps)}
-
q^{(kp+r)(kp'+s)}\right)\\
&=&
\sum
\left(
q^{n_{1,d} + \sum_{i=1}^{f_1} n_{1,i}^2
+\sum_{k=2}^{x} \sum_{i=1}^{f_k} n_{k,i}(n_{k,i}-1)}
\phantom{M\atop N}\right. \\[0.5mm]
&&\left.\qquad\times\phantom{M\atop N}
{1\over (q)_{2n_{2,f_2}}}
\prod_{i=2}^{f_1} {1\over (q)_{n_{1,i-1}-n_{1,i}}}
\prod_{i=1}^{f_2} {1\over (q)_{n_{2,i-1}-n_{2,i}}}
\right. \\[0.5mm]
&&\left.\qquad\times\phantom{M\atop N}
\prod_{k=3}^{x}\prod_{i=1+\delta_{3,k}}^{f_k}
\left[{2\left(w_k-\sum_{j=1}^i n_{k,j}\right)+n_{k,i-1}+n_{k,i}
\atop n_{k,i-1}-n_{k,i}}\right]
\right),
\end{eqnarray*}
where the sums are taken over the same parameters as in Theorem \ref{DinfThrm}.
\end{corollary}

\section{Discussion}

$Q$-polynomial identities for the Virasoro characters that 
we are interested in, were previously studied in \cite{bm,bms}. 
We wish to discuss the similarities and differences between 
this work and \cite{bm,bms}.

In this work, we have used purely combinatorial methods: the Burge 
transform, plus a combinatorial identity. This enabled us to  
obtain two polynomial identities for each character 
$\chi^{p, p'}_{r, s}$ in the subset discussed above. Obtaining
identities for more characters would require either an
extension of the Burge transform, or further combinatorial 
identities. 

In \cite{bm,bms}, Bethe Ansatz type methods were used, and 
an identity for each character in various subsets of 
$\chi^{p, p'}_{r, s}$ was obtained. It is our understanding 
that the methods of \cite{bm,bms} are sufficiently general 
to produce polynomial identities for {\it all} characters.
However, this task is impeded by computational complexity. 

Aside from the fact that our polynomial identities that depend 
on two finite size parameters have no counterparts in
\cite{bm,bms}, there are overlaps between our results and those 
of \cite{bm,bms} at the level of identities with one finite size 
parameter.

Our results of the latter type are of two forms: identities 
with a conventional bosonic side, and identities with a
non-conventional bosonic side. The results of \cite{bm,bms}
are all of the former type. In this case, our results form 
only a subset of those of \cite{bm,bms}.

For example, it is easy to show that we obtain the same identities 
as in \cite{bm,bms} for $\chi^{p, p'}_{r_0, s_0}$, where 
$| p s_0 - p' r_0 | = 1$. On the other hand, all identities that we 
obtain and that have a non-conventional bosonic side are new. 

\bigskip

{\footnotesize
{\noindent \it Acknowledgements --- }

We are indebted to G. E. Andrews, W. H. Burge, and Ch. Krattenthaler 
for drawing our attention to \cite{burge}, and for correspondence and 
clarifications. We also wish to thank A. Berkovich, B. M. McCoy, and  
A. Schilling for discussions and explanations of \cite{bms}. This work 
was made possible by financial support of the Australian Research Council 
(ARC).}

\begin{appendix}

\section{Derivation}

In this appendix, we derive Theorem \ref{GenFunThrm}.
To do this, we classify each element of
${\cal R}(s,u,t,v,a,b,\alpha,\beta)$ according to its
{\it sequence of faults} as follows.
First note that the conditions (\ref{restrictions}) and
(\ref{exceptions}) for a partition pair
$(q,p)$ in $(s,u)\times(t,v)$ to be restricted may be expressed:
\begin{eqnarray*}
p_i-q_{i+1-a} &\ge& 1-\alpha\qquad(a\le i\le t)\cr
q_i-p_{i+1-b} &\ge& 1-\beta\qquad(b\le i\le s),
\end{eqnarray*}
where, when required, we take $q_0=p_0=0$, $q_{s+1}=u$
and $p_{t+1}=v$.
For $(q,p)\in{\cal R}(s,u,t,v,a,b,\alpha,\beta)$,
locate the largest $i$ that violates these conditions.
Note that if $a>0$ and $b>0$ then both conditions cannot be
violated simultaneously because
if the first is violated so that $q_{i+1-a}>p_i+\alpha$ then for
$j\ge i+1-a$ and $k\le i$,
$$
q_j\ge q_{i+1-a}\ge p_i+\alpha \ge p_k+\alpha \ge p_k-\beta+1
$$
since $\alpha+\beta\ge1$, thereby, in the particular case
where $j=i$ and $k=i+1-b$,
satisfying the second condition above.
The violation at position $i$ is said to be an $a$-fault if it is the first
of the two conditions above that is violated. 
Otherwise it is said to be a $b$-fault.
Now locate the next largest $i$ for which a fault of the {\it opposite}
variety occurs. (Note that if the first fault is an $a$-fault at position
$i_1$ then the second fault will be at position $i_2$ satisfying
$i_2\le i_1-a$, and will be a $b$-fault. On the other hand if the
first fault is a $b$ fault at position $i_1$ then the second will
be an $a$-fault at position $i_2\le i_1-b$.)
In this way, produce an alternating sequence of $a$s and $b$s for
each element of ${\cal R}(s,u,t,v,a,b,\alpha,\beta)$.

If either $a=0$ or $b=0$ then it is possible for an $a$-fault and a
$b$-fault to occur at the same $i$. In this case, both faults
must be recorded. If $a=0$, we consider the $a$-fault as preceding
the $b$-fault, whereas if $b=0$, we consider the $b$-fault as
preceding the $a$-fault.

Define ${\cal A}_k(s,u,t,v,a,b,\alpha,\beta)$ to be the set that
comprises all those
elements of ${\cal R}(s,u,t,v,a,b,\alpha,\beta)$ whose sequence of
faults {\it contains} a subsequence of length $k$ starting with an $a$.
Likewise, define
the set ${\cal B}_k(s,u,t,v,a,b,\alpha,\beta)$ to comprise all those
elements of ${\cal R}(s,u,t,v,a,b,\alpha,\beta)$ whose sequence of
faults {\it contains} a subsequence of length $k$ starting with a $b$.
Note that for each $k$, ${\cal A}_k \subset {\cal A}_{k-1}$ and
${\cal B}_k \subset {\cal B}_{k-1}$.
The generating functions for these sets are defined by:
\begin{eqnarray*}
A_k(s,u,t,v,a,b,\alpha,\beta)&=&
\sum_{\pi\in{\cal A}_k(s,u,t,v,a,b,\alpha,\beta)}
q^{\wt(\pi)},\cr
B_k(s,u,t,v,a,b,\alpha,\beta)&=&
\sum_{\pi\in{\cal B}_k(s,u,t,v,a,b,\alpha,\beta)}
q^{\wt(\pi)}.
\end{eqnarray*}

\begin{lemma}\label{StepLemma}
\begin{eqnarray*}
A_k(s,u,t,v,a,b,\alpha,\beta)&=&
q^{\alpha(s-t+a)}
B_{k-1}(t-a,v+\alpha,s+a,u-\alpha,a,b,\alpha,\beta),\cr
B_k(s,u,t,v,a,b,\alpha,\beta)&=&
q^{\beta(t-s+b)}
A_{k-1}(t+b,v-\beta,s-b,u+\beta,a,b,\alpha,\beta).
\end{eqnarray*}
\end{lemma}

\Proof The first of these results is proved by setting up a
bijection between the sets
${\cal A}_k(s,u,t,v,a,b,\alpha,\beta)$ and
${\cal B}_{k-1}(t-a,v+\alpha,s+a,u-\alpha,a,b,\alpha,\beta)$.
For each $(q,p)\in{\cal A}_k(s,u,t,v,a,b,\alpha,\beta)$, the
bijective image is obtained as follows.
Let $i$ be the position of the first $a$-fault of $(q,p)$,
so that $p_i-q_{i-a+1}<1-\alpha$, whereupon $p_i\le q_{i-a+1}-\alpha$.
In addition, since this is the first $a$-fault, we have
$p_{i+1}-q_{i-a+2}\ge 1-\alpha$ whereupon
$p_{i+1}+\alpha>q_{i-a+2}\ge q_{i-a}$.
(Actually, care should be taken at the exceptional values. First note
that if $t-s=a$ then $p_t\ge u-\alpha+1=p_{s+1}-\alpha+1$ implies that
$p_t-q_{t-a+1}\ge1-\alpha$ so that no $a$-fault then occurs at $i=t$.
Thus either $i<t$ or $t-s<a$. Thus if $i=t$ then, from $u-v\le\alpha$
follows $q_{i-a}=q_{t-a}\le q_s\le u\le v+\alpha$.
Also if $i=t-1$ and $t-s=a$ then
$q_{i-a}=q_{t-a-1}\le q_s\le u\le p_t+\alpha-1<p_{i+1}+\alpha$.
Otherwise $i+1\le t$ and $i-a+2\le s$ and the previous reasoning holds.
Also note that if $a=0$ then $q_1\le\alpha-1$ implies that
$p_0-q_{-a+1}\ge1-\alpha$ so that no $a$-fault occurs at $i=0$.
Thus $i\ge\max\{1,a\}$.)
These two inequalities ensure that $p^\prime$ and $q^\prime$
defined by:
\bea
p^\prime&=&(q_s-\alpha,q_{s-1}-\alpha,\ldots,
q_{i+1-a}-\alpha,p_i,p_{i-1},\ldots,p_1),\\[0.5mm]
q^\prime&=&(p_t+\alpha,p_{t-1}+\alpha,\ldots,
p_{i+1}+\alpha,q_{i-a},q_{i-a-1},\ldots,q_1)
\eea
are both partitions. Furthermore, the partition pair $(q^\prime,p^\prime)$
is in $(t-a,v+\alpha)\times(s+a,u-\alpha)$.
Now the faults beyond the $i$th position in $(q^\prime,p^\prime)$
are precisely as they are in $(q,p)$.
This ensures that $(q^\prime,p^\prime)\in
{\cal B}_{k-1}(t-a,v+\alpha,s+a,u-\alpha,a,b,\alpha,\beta)$.

The reverse map is defined as follows.
Let $(q^\prime,p^\prime)\in
{\cal B}_{k-1}(t-a,v+\alpha,s+a,u-\alpha,a,b,\alpha,\beta)$
and let $i$ be the position of the first $b$-fault so that
$q^\prime_i-p^\prime_{i-b+1}<1-\beta$.
Now locate $j\ge i+a$ such that $p^\prime_{j+1}+\alpha\ge q^\prime_{j-a}$ and
$p^\prime_{h+1}+\alpha<q^\prime_{h-a}$ for $h>j$.
If such a $j$ exists then the second inequality implies that
$p^\prime_j\le p^\prime_{j+2}<q^\prime_{j-a+1}-\alpha$,
thereby ensuring that $p$ and $q$ defined by
\bea
p&=&(q^\prime_{t-a}-\alpha,q^\prime_{t-a-1}-\alpha,\ldots,
q^\prime_{j-a+1}-\alpha,p^\prime_j,p^\prime_{j-1},\ldots,p^\prime_1),\\[0.5mm]
q&=&(p^\prime_{s+a}+\alpha,p^\prime_{s+a-1}+\alpha,\ldots,
p^\prime_{j+1}+\alpha,q^\prime_{j-a},q^\prime_{j-a-1},\ldots,q^\prime_1)
\eea
are each partitions.
Furthermore $(q,p)$ is in $(s,u)\times(t,v)$.
In addition to the faults at and to the right of the $i$th position,
$(q,p)$ has an $a$-fault at position $j\ge i+a$ since
$p_{j}-q_{j-a+1}=p^\prime_j-p^\prime_{j+1}-\alpha\le-\alpha<1-\alpha$.
Therefore, in the case where $j$ can be found, the partition pair
$(q,p)$ is an element of ${\cal A}_k(s,u,t,v,a,b,\alpha,\beta)$.
Also note that if $k>j$ then
$p_k-q_{k-a+1}=q^\prime_{k-a}-\alpha-(p^\prime+\alpha)>
q^\prime_{k-a}-\alpha-q^\prime_{k-a}\ge 1-\alpha$, thus ensuring
that the $a$-fault of $(q,p)$ at position $j$ is the first,
and therefore that the map described in this paragraph is
the inverse of that given in the previous paragraph.

In the case where such a $j$ cannot be located, so that
$p^\prime_{j+1}+\alpha<q^\prime_{j-a}$ for all $j\ge i+a$,
define
\bea
p&=&(p^\prime_t,p^\prime_{t-1},p^\prime_{t-2},\ldots,p^\prime_1),\\[0.5mm]
q&=&(p^\prime_{s+a}+\alpha,p^\prime_{s+a-1}+\alpha,\ldots,
p^\prime_{t+1}+\alpha,q^\prime_{t-a},q^\prime_{t-a-1},\ldots,q^\prime_1),
\eea
whence an $a$-fault has been introduced at position $t$ since
$p_t-q_{t-a+1}=p^\prime_t-p^\prime_{t+1}-\alpha\le -\alpha<1-\alpha$.
The faults of $(q^\prime,p^\prime)$ also occur in $(q,p)$
to the right of position $s-a$.
Thus also in this case
$(q,p)\in{\cal A}_k(s,u,t,v,a,b,\alpha,\beta)$.
Certainly, the introduced $a$-fault is the first fault so that in this
case, this map is the inverse of that above.

Thus the two sets ${\cal A}_k(s,u,t,v,a,b,\alpha,\beta)$ and
${\cal B}_{k-1}(t-a,v+\alpha,s+a,u-\alpha,a,b,\alpha,\beta)$
are in bijection.
The first expression of the lemma then follows after noting that if
$(q^\prime,p^\prime)$ is the bijective image of $(q,p)$
then $\wt((q,p))=\wt((q^\prime,p^\prime))+(s-t+a)\alpha$.

The second expression is proved in a totally analogous manner.
\cqfd

\begin{lemma}\label{IndLem}
\begin{eqnarray*}
\hbox to 5mm{$A_{2k}(s,u,t,v,a,b,\alpha,\beta)$\hss}\\[0.4mm]
&=&q^{k^2(a+b)(\alpha+\beta)+k(\alpha+\beta)(s-t)+k(a\beta-b\alpha)}
g(s,u,t,v,k(a+b),k(\alpha+\beta))\\[0.5mm]
\hbox to 5mm{$B_{2k}(s,u,t,v,a,b,\alpha,\beta)$\hss}\\[0.4mm]
&=&q^{k^2(a+b)(\alpha+\beta)-k(\alpha+\beta)(s-t)-k(a\beta-b\alpha)}
g(t,v,s,u,k(a+b),k(\alpha+\beta))\\[0.5mm]
\hbox to 5mm{$A_{2k+1}(s,u,t,v,a,b,\alpha,\beta)$\hss}\\[0.4mm]
&=&q^{k^2(a+b)(\alpha+\beta)+(k(\alpha+\beta)+\alpha)(s-t)
+k(a\beta+b\alpha+2a\alpha)+a\alpha}\\[0.4mm]
&&\qquad\qquad\qquad\qquad\times\quad
g(s,u,t,v,k(a+b)+a,k(\alpha+\beta)+\alpha)\\[0.5mm]
\hbox to 5mm{$B_{2k+1}(s,u,t,v,a,b,\alpha,\beta)$\hss}\\[0.4mm]
&=&q^{k^2(a+b)(\alpha+\beta)-(k(\alpha+\beta)+\beta)(s-t)
+k(a\beta+b\alpha+2b\beta)+b\beta}\\[0.4mm]
&&\qquad\qquad\qquad\qquad\times\quad
g(t,v,s,u,k(a+b)+b,k(\alpha+\beta)+\beta).
\end{eqnarray*}
\end{lemma}

\Proof
We proceed by induction. The first two results clearly hold when
$k=0$ since the required generating function is that
for all partitions, regardless of their faults.
Therefore
\begin{eqnarray*}
\hbox to 5mm{$A_0(s,u,t,v,a,b,\alpha,\beta)
=B_0(s,u,t,v,a,b,\alpha,\beta)$\hss}\\[0.5mm]
&&\qquad=\left[s+u\atop s\right]\left[t+v\atop t\right]
=g(s,u,t,v,0,0)=g(t,v,s,u,0,0).
\end{eqnarray*}

Now assume that for a given $i$, the expressions hold for
$A_i(s,u,t,v,a,b,\alpha,\beta)$ and $B_i(s,u,t,v,a,b,\alpha,\beta)$.
We show that this implies that the expressions for
$A_{i+1}(s,u,t,v,a,b,\alpha,\beta)$ and $B_{i+1}(s,u,t,v,a,b,\alpha,\beta)$
hold.
In the case where $i$ is even, let $i=2k$ whereupon, on using
Lemma \ref{StepLemma}, we obtain:
\begin{eqnarray*}
\hbox to 5mm{$A_{i+1}(s,u,t,v,a,b,\alpha,\beta)$\hss}\\[0.5mm]
&=&
q^{\alpha(s-t+a)}B_i(t-a,v+\alpha,s+a,u-\alpha,a,b,\alpha,\beta)\\[0.5mm]
&=&
q^{\alpha(s-t+a)
+k^2(a+b)(\alpha+\beta)-k(\alpha+\beta)(t-s-2a)-k(a\beta-b\alpha)}\\[0.5mm]
&&\qquad\qquad\qquad\times\quad
g(s+a,u-\alpha,t-a,v+\alpha,k(a+b),k(\alpha+\beta))\\[0.5mm]
&=&
q^{a\alpha+k^2(a+b)(\alpha+\beta)+(k(\alpha+\beta)+\alpha)(s-t)
+k(a\beta+b\alpha+2a\alpha)}\\[0.5mm]
&&\qquad\qquad\qquad\times\quad
g(s,u,t,v,k(a+b)+a,k(\alpha+\beta)+\alpha),
\end{eqnarray*}
where use has been made of
\begin{eqnarray*}
\hbox to 5mm{$g(s+a,u-\alpha,t-a,v+\alpha,x,y)$\hss}\\[0.5mm]
&=&\left[{s+a+u-\alpha+x-y\atop s+a+x}\right]
\left[{t-a+v+\alpha-x+y\atop t-a-x}\right]\\[0.5mm]
&=&g(s,u,t,v,x+a,y+\alpha).
\end{eqnarray*}
Similarly,
\begin{eqnarray*}
\hbox to 5mm{$B_{i+1}(s,u,t,v,a,b,\alpha,\beta)$\hss}\\[0.5mm]
&=&
q^{\beta(t-s+b)}A_i(t+b,v-\beta,s-b,u+\beta,a,b,\alpha,\beta)\\[0.5mm]
&=&
q^{\beta(t-s+b)
+k^2(a+b)(\alpha+\beta)+k(\alpha+\beta)(t-s+2b)+k(a\beta-b\alpha)}\\[0.5mm]
&&\qquad\qquad\qquad\times\quad
g(t+b,v-\beta,s-b,u+\beta,k(a+b),k(\alpha+\beta))\\[0.5mm]
&=&
q^{b\beta+k^2(a+b)(\alpha+\beta)-(k(\alpha+\beta)+\beta)(s-t)
+k(a\beta+b\alpha+2b\beta)}\\[0.5mm]
&&\qquad\qquad\qquad\times\quad
g(t,v,s,u,k(a+b)+b,k(\alpha+\beta)+\beta).
\end{eqnarray*}
In the case where $i$ is odd, let $i=2k+1$ whereupon,
on using Lemma \ref{StepLemma},
\begin{eqnarray*}
\hbox to 5mm{$A_{i+1}(s,u,t,v,a,b,\alpha,\beta)$\hss}\\[0.5mm]
&=&
q^{\alpha(s-t+a)}B_i(t-a,v+\alpha,s+a,u-\alpha,a,b,\alpha,\beta)\\[0.5mm]
&=&
q^{\alpha(s-t+a)
+k^2(a+b)(\alpha+\beta)-(k(\alpha+\beta)+\beta)(t-s-2a)
+k(a\beta+b\alpha+2b\beta)+b\beta}\\[0.5mm]
&&\qquad\qquad\times\quad
g(s+a,u-\alpha,t-a,v+\alpha,k(a+b)+b,k(\alpha+\beta)+\beta)\\[0.5mm]
&=&
q^{k^2(a+b)(\alpha+\beta)+(k+1)(\alpha+\beta)(s-t)
+2ka(\alpha+\beta)+k(a\beta+b\alpha+2b\beta)+2a\beta+a\alpha+b\beta}\\[0.5mm]
&&\qquad\qquad\times\quad
g(s,u,t,v,(k+1)(a+b),(k+1)(\alpha+\beta))\\[0.5mm]
&=&
q^{(k+1)^2(a+b)(\alpha+\beta)+(k+1)(\alpha+\beta)(s-t)
+(k+1)(a\beta-b\alpha)}\\[0.5mm]
&&\qquad\qquad\times\quad
g(s,u,t,v,(k+1)(a+b),(k+1)(\alpha+\beta)).
\end{eqnarray*}
Similarly,
\begin{eqnarray*}
\hbox to 5mm{$B_{i+1}(s,u,t,v,a,b,\alpha,\beta)$\hss}\\[0.5mm]
&=&
q^{\beta(t-s+b)}A_i(t+b,v-\beta,s-b,u+\beta,a,b,\alpha,\beta)\\[0.5mm]
&=&
q^{\beta(t-s+b)
+k^2(a+b)(\alpha+\beta)+(k(\alpha+\beta)+\alpha)(t-s+2b)
+k(a\beta+b\alpha+2a\alpha)+a\alpha}\\[0.5mm]
&&\qquad\qquad\times\quad
g(t+b,v-\beta,s-b,u+\beta,k(a+b)+a,k(\alpha+\beta)+\alpha)\\[0.5mm]
&=&
q^{k^2(a+b)(\alpha+\beta)-(k+1)(\alpha+\beta)(s-t)
+2kb(\alpha+\beta)+k(a\beta+b\alpha+2a\alpha)+2b\alpha+a\alpha+b\beta}\\[0.5mm]
&&\qquad\qquad\times\quad
g(t,v,s,u,(k+1)(a+b),(k+1)(\alpha+\beta))\\[0.5mm]
&=&
q^{(k+1)^2(a+b)(\alpha+\beta)-(k+1)(\alpha+\beta)(s-t)
-(k+1)(a\beta-b\alpha)}\\[0.5mm]
&&\qquad\qquad\times\quad
g(t,v,s,u,(k+1)(a+b),(k+1)(\alpha+\beta)).
\end{eqnarray*}
Thus for both even and odd $i$, the expressions for
$A_{i+1}$ and $B_{i+1}$ follow from those for $A_{i}$ and $B_{i}$.
The lemma then follows by induction.
\cqfd

We are now in a position to give a {\it non-constant sign} expression
for the generating function $R(s,u,t,v,a,b,\alpha,\beta)$.

\begin{theorem}\label{BosGenThrm}
\begin{eqnarray*}
\hbox to 5mm{$R(s,u,t,v,a,b,\alpha,\beta)$\hss}\\[1mm]
&=& \sum_{k=-\infty}^\infty
q^{k^2(a+b)(\alpha+\beta)+k(\alpha+\beta)(s-t)+k(a\beta-b\alpha)}\\[-4mm]
&&\qquad\qquad\qquad\qquad\quad\times\quad
g(s,u,t,v,k(a+b),k(\alpha+\beta))\\[1mm]
&&\quad-\quad \sum_{k=-\infty}^\infty
q^{k^2(a+b)(\alpha+\beta)+(k(\alpha+\beta)+\alpha)(s-t)
+k(a\beta+b\alpha+2a\alpha)+a\alpha}\\[-4mm]
&&\qquad\qquad\qquad\qquad\quad\times\quad
g(s,u,t,v,k(a+b)+a,k(\alpha+\beta)+\alpha).
\end{eqnarray*}
\end{theorem}

\Proof
We require the generating function for partition pairs
in $(s,u)\times(t,v)$ that have no faults.
In this proof for typographical reasons, we drop the
arguments $(s,u,t,v,a,b,\alpha,\beta)$
from $R$, $A_i$, $B_i$, ${\cal R}$, ${\cal A}_i$ and ${\cal B}_i$.
Note first that ${\cal A}_{2k}\subset{\cal A}_{2k-1}$.
Moreover, the set ${\cal A}_{2k-1}\backslash{\cal A}_{2k}$
contains all partition pairs that have a sequence of
faults $abab\cdots ba$ of length $2k-1$ or a sequence
of faults $baba\cdots ba$ of length $2k$.
Likewise, the set ${\cal B}_{2k-1}\backslash{\cal B}_{2k}$
contains all partition pairs that have a sequence of
faults $baba\cdots ab$ of length $2k-1$ or a sequence
of faults $abab\cdots ab$ of length $2k$.

Therefore, $({\cal A}_{2k-1}\backslash{\cal A}_{2k})\cup
({\cal B}_{2k-1}\backslash{\cal B}_{2k})$ contains all partition
pairs with a sequence of faults of length $2k$ or $2k-1$.
The generating function for such sequences is therefore
$A_{2k-1}-A_{2k}+B_{2k-1}-B_{2k}$.
Thereupon, the generating function for sequences with no
faults is:
$$
R=A_0-\sum_{k=1}^\infty (A_{2k-1}-A_{2k}+B_{2k-1}-B_{2k}).
$$
The theorem then follows from the four expressions given by
Lemma \ref{IndLem} after noting that
$$
g(t,v,s,u,x,y)
=\left[{t+v+x-y\atop t+x}\right] \left[{s+u-x+y\atop s-x}\right]
=g(s,u,t,v,-x,-y),
$$
and that the fourth expression of Lemma \ref{IndLem} may be
re-expressed
\begin{eqnarray*}
\hbox to 5mm{$B_{2k-1}(s,u,t,v,a,b,\alpha,\beta)$\hss}\\[1mm]
&=&q^{k^2(a+b)(\alpha+\beta)-(k(\alpha+\beta)-\alpha)(s-t)
-k(a\beta+b\alpha+2a\alpha)+a\alpha}\\[0.5mm]
&&\qquad\qquad\qquad\qquad\quad\times\quad
g(s,u,t,v,-k(a+b)+a,-k(\alpha+\beta)+\alpha).
\end{eqnarray*}
\cqfd

\section{Details of Examples}\label{MinfApp}

\begin{example}\label{MinfExAp} {\rm
Here, we provide the details for Example \ref{MinfEx}.
The following table lists the partition sequences used to calculate
$\lim_{M\to\infty}R(3,M,2,14,5,32)$ by means of Theorem \ref{MinfThrm}.
\medskip
$$
\vbox{\offinterlineskip\halign{$\strut
\enspace #\hfil\enspace$&\vrule #&&$\hfil\enspace #\hfil$\cr
(n^{(1)};n^{(2)};\cdots;n^{(x)})
&&x&n_{1,0}&w_1&n_{2,0}&w_2&n_{3,0}&w_3&n_{4,0}&w_4\cr
\omit&height 1mm\cr
\multispan{11}\hrulefill\cr
\omit&height 1mm\cr
(3,3;3,0)&&2&3&\infty&\infty&3\cr
(3,2;2,0)&&2&3&\infty&\infty&2\cr
(3,1;1,0)^*&&2&3&\infty&\infty&1\cr
(3,0;0,0)^*&&2&3&\infty&\infty&0\cr
(2,2;2,0)^*&&2&3&\infty&\infty&2\cr
(2,1;1,0)^*&&2&3&\infty&\infty&1\cr
(2,0;0,0)^*&&2&3&\infty&\infty&0\cr
(1,1;1,0)^*&&2&3&\infty&\infty&1\cr
(1,0;0,0)^*&&2&3&\infty&\infty&0\cr
(0,0;0,0)^*&&2&3&\infty&\infty&0\cr
(3,3;2,0;0)&&3&3&\infty&\infty&3&1&0\cr
(3,3;1,1;0)&&3&3&\infty&\infty&3&2&0\cr
(3,2;1,0;0)&&3&3&\infty&\infty&2&1&0\cr
(3,1;0,0;0)^*&&3&3&\infty&\infty&1&1&0\cr
(2,2;1,0;0)^*&&3&3&\infty&\infty&2&1&0\cr
(2,1;0,0;0)^*&&3&3&\infty&\infty&1&1&0\cr
(1,1;0,0;0)^*&&3&3&\infty&\infty&1&1&0\cr
(3,3;1,0;1;0,0,0)&&4&3&\infty&\infty&3&1&1&1&0\cr
(3,3;0,0;1;0,0,0)&&4&3&\infty&\infty&3&1&2&2&0\cr
(3,2;0,0;1;0,0,0)&&4&3&\infty&\infty&2&1&1&1&0\cr
(2,2;0,0;1;0,0,0)^*&&4&3&\infty&\infty&2&1&1&1&0\cr
}}
$$
\medskip
\noindent
(Those sequences marked with an asterisk are also relevant to Example
\ref{DinfEx}.)
The use of Theorem \ref{MinfThrm} then gives:
\begin{eqnarray*}
\hbox to 0mm{$\displaystyle
\lim_{M\to\infty}R(3,M,2,14,5,32)
$\hss}\\[0.5mm]
&=&
q^{27} \qeta6 \qeta0 \qeta0 \gauss33
+ q^{17} \qeta4 \qeta0 \qeta1 \gauss22
+ q^{11} \qeta2 \qeta0 \qeta2 \gauss11
\\[0.5mm]
&&
\:{} + q^9    \qeta0 \qeta0 \qeta3 \gauss00
+ q^{12} \qeta4 \qeta1 \qeta0 \gauss22
+ q^6    \qeta2 \qeta1 \qeta1 \gauss11
\\[0.5mm]
&&
\:{} + q^4    \qeta0 \qeta1 \qeta2 \gauss00
+ q^3    \qeta2 \qeta2 \qeta0 \gauss11
+ q      \qeta0 \qeta2 \qeta1 \gauss00
\\[0.5mm]
&&
\:{} +        \qeta0 \qeta3 \qeta0 \gauss00
+ q^{22} \qeta6 \qeta0 \qeta0 \gauss42 \gauss11
+ q^{20} \qeta6 \qeta0 \qeta0 \gauss40 \gauss22
\\[0.5mm]
&&
\:{} + q^{14} \qeta4 \qeta0 \qeta1 \gauss31 \gauss11
+ q^{10} \qeta2 \qeta0 \qeta2 \gauss20 \gauss11
\\[0.5mm]
&&
\:{} + q^9    \qeta4 \qeta1 \qeta0 \gauss31 \gauss11
+ q^5    \qeta2 \qeta1 \qeta1 \gauss20 \gauss11
+ q^2    \qeta2 \qeta2 \qeta0 \gauss20 \gauss11
\\[0.5mm]
&&
\:{} + q^{19} \qeta6 \qeta0 \qeta0 \gauss51 \gauss20 \gauss11 \gauss00 \gauss00
+ q^{18} \qeta6 \qeta0 \qeta0 \gauss60 \gauss40 \gauss22 \gauss00 \gauss00
\\[0.5mm]
&&
\:{} + q^{13} \qeta4 \qeta0 \qeta1 \gauss40 \gauss20 \gauss11 \gauss00 \gauss00
+ q^8    \qeta4 \qeta1 \qeta0 \gauss40 \gauss20 \gauss11 \gauss00 \gauss00
\\[0.5mm]
&=&
1+2q+5q^2+10q^3+18q^4+30q^5+49q^6+74q^7+110q^8+158q^9+\cdots
\end{eqnarray*}
\endex}
\end{example}

\end{appendix}


\begin{thebibliography}{99}

\bibitem{bpz}{\sc A.~A.~Belavin, A.~M.~Polyakov} and
             {\sc A.~B.~Zamolodchikov}, 
             {\it Infinite conformal symmetry in 
	          two-dimensional quantum field theory}, 
              Nucl. Phys. {\bf B 241} (1984) 333. 

\bibitem{dms-book}{\sc Ph.~Di~Francesco, P.~Mathieu} and 
                {\sc D.~Senechal},
	        {\it Conformal Field Theory},
	             Springer, 1996.

\bibitem{abf}{\sc G.E.~Andrews, R.J.~Baxter} and {\sc P.J.~Forrester},
{\it Eight-vertex SOS model and generalized Rogers-Ramanujan-type 
     identities},
     J. Stat. Phys. {\bf 35} (1984), 193.

\bibitem{fb}{\sc P.J.~Forrester} and {\sc R.J.~Baxter},
{\it Further exact solutions of the eight-vertex SOS model and
generalizations of the Rogers-Ramanujan identities},
J. Stat. Phys. {\bf 38} (1985) 435. 

\bibitem{baxter-book}{\sc R.~J.~Baxter},
             {\it Exactly solvable models in statistical mechanics}, 
	     Academic Press, London, 1982.

\bibitem{djkmo}{\sc E.~Date, M.~Jimbo, A.~Kuniba, T.~Miwa} and 
               {\sc M.~Okado},
               {\em Exactly solvable SOS models: Local height 
	            probabilities and theta function identities},
                    Nucl. Phys. {\bf B 290 [FS20]} (1987) 231.

\bibitem{stony-brook-review}{\sc R.~Kedem, B.~M.~McCoy} and {\sc E.~Melzer}, 
{\it The sums of Rogers, Schur and Ramanujan and the Bose-Fermi
correspondence in $1+1$ dimensional quantum field theory}, 
{\tt hep-th 9304056}.

\bibitem{lepowsky}{\sc J.~Lepowsky and M.~Primc},
                  {\it Structure of the standard modules of the 
		       affine Lie algebra A$_1^{(1)}$},
                  {\it Contemporary Mathematics}, {\bf 46} 
		       (AMS, Providence, 1985).

\bibitem{rocha}{\sc A.~Rocha-Caridi},
        {\em Vacuum vector representations of the Virasoro
	algebra}, in
	Vertex Operators in Mathematics and Physics,
	eds. J.~Lepowsky, S.~Mandelstam and I.M.~Singer,
	Springer, 1985.

\bibitem{andrews-red-book}{\sc G.E.~Andrews},
                          {\it The Theory of Partitions},
                               Encyclopedia of Mathematics and 
			       its Applications, 
	                  {\bf 2},
                               Addison-Wesley, 1976.

\bibitem{abbbfv}{\sc G.E.~Andrews, 
                     R.J.~Baxter, 
                     D.M.~Bressoud, 
		     W.H.~Burge, 
		     P.J.~Forrester}, and 
                {\sc G.X.~Viennot},
                {\it Partitions with prescribed hook differences},
                     European J. of Combin. {\bf 8} (1987) 4310.

\bibitem{burge}{\sc W.~H.~Burge},
               {\it Restricted partition pairs},
	       {\it J. of Comb. Th.} {\bf A~63}
	       (1993) 210.

\bibitem{bailey}{\sc W.~N.~Bailey}, 
                {\it Proc. London Math. Soc.} {\bf 50} (1949) 1.

\bibitem{andrews-on-bailey} {\sc G.~E.~Andrews},
                            {\it Multiple series Rogers-Ramanujan 
			         type identities},
                            Pacific J. Math. {\bf 114} (1984) 267.

\bibitem{andrews-white-book}{\sc G.~E.~Andrews},
        {\em $q$-Series: Their development and application in 
	     analysis, number theory, combinatorics, physics, 
	     and computer algebra},
             in CBMS Regional Conf. Ser. in Math., {\bf 66},
            (AMS, Providence, Rhode Island, 1985).

\bibitem{fq} {\sc O.~Foda} and {\sc Y.-H.~Quano},
             {\it Virasoro character identities from 
	          the Andrews--Bailey construction}, 
                  Int.\ J. Mod.\ Phys.\ A {\bf 12} (1997) 1651.

\bibitem{bms}{\sc A.~Berkovich, B.~M.~McCoy and A.~Schilling},
             {\it Rogers--Schur--Ramanujan type identities for 
	          the $M(p, p')$ minimal models of conformal 
		  field theory}, 
                  preprint {\tt q-alg/9607020},
		  to appear in Commun.\  Math.\  Phys.

\bibitem{bm}{\sc A.~Berkovich and B.~M.~McCoy},
            {\it Continued fractions and fermionic representations 
	         for characters of $M(p, p')$ minimal models}, 
                 Lett. \ Math.\ Phys. {\bf 37} (1996) 49.

\bibitem{fw}{\sc O.~Foda} and {\sc S.O.~Warnaar},
            {\it A bijection which implies Melzer\rq s polynomial
	         identities: the $\chi^{(p,p+1)}_{1,1}$ case},
	         Lett.\ Math.\ Phys.\ {\bf 36} (1996) 145.

\bibitem{w1}{\sc S.O.~Warnaar},
            {\it Fermionic solution of the Andrews-Baxter-Forrester
	         model I. Unification of TBA and CTM methods},
                 J.\ Stat.\ Phys.\ {\bf 82} (1996) 657.

\bibitem{w2}{\sc S.O.~Warnaar},
            {\it Fermionic solution of the Andrews-Baxter-Forrester
	         model. II. Proof of Melzer\rq s polynomial
		 identities},
	         J.\ Stat.\ Phys.\ {\bf 84} (1997) 49.

\bibitem{fow}{\sc O.~Foda, M.~Okado} and {\sc S.O.~Warnaar},
        {\em A proof of polynomial identities of type
	     $\widehat{sl(n)_1}\otimes\widehat{sl(n)_1}
	     /\widehat{sl(n)_2}$},
             J. Math. Phys. {\bf 37} (1996) 965.

\bibitem{kkr} {\sc S.~Kerov, A.N.~Kirillov} and {\sc
                 N.~Reshetikhin}, J. Sov. Math. {\bf 41} (1988)
		 916.

\bibitem{df}{\sc S.~Dasmahapatra} and {\sc O.~Foda},
            {\it Strings, paths and standard tableaux},
            {\tt q-alg 9601011}, 
                 to appear in Int.\ J.\ of Mod.\ Phys.\ (1997).

\bibitem{kg} {\sc Ch.~Krattenthaler} and {\sc I.~Gessel}, 
             {\it Cylindric partitions}, 
	          Trans.\  Amer.\ Math.\ Soc. {\bf 349} (1997), 429-479.

\bibitem{gould}{\sc H.W.~Gould},
        {\em A new symmetrical combinatorial identity},
	J. Comb. Theory (A) {\bf 13} (1972) 278.

\end{thebibliography}
\end{document}